\documentclass[a4paper,fleqn]{cas-sc}

\usepackage{setspace}
\usepackage{bm}
\usepackage{multirow}
\usepackage[numbers,sort&compress]{natbib}
\usepackage{graphicx}
\usepackage{float}
\usepackage{geometry}
\usepackage{enumitem}
\usepackage{amsmath}
\usepackage{times}
\usepackage{caption}
\usepackage[ruled]{algorithm2e}
\usepackage[title]{appendix}%

\usepackage{lineno}
\def\tsc#1{\csdef{#1}{\textsc{\lowercase{#1}}\xspace}}
\tsc{WGM}
\tsc{QE}
\tsc{EP}
\tsc{PMS}
\tsc{BEC}
\tsc{DE}
\begin{document}

\let\WriteBookmarks\relax
\def\floatpagepagefraction{1}
\def\textpagefraction{.001}
\shorttitle{Analytical Pursuit-Evasion Game Strategy in Arbitrary Keplerian Reference Orbits}
\shortauthors{Shuyue Fu et~al.}

\title [mode = title]{Analytical Pursuit-Evasion Game Strategy in Arbitrary Keplerian Reference Orbits}                      

\author[1,2]{Shuyue Fu}[orcid=0009-0001-5111-9779,style=chinese]
\ead{fushuyue@buaa.edu.cn}
\credit{Data curation, Formal analysis, Methodology, Software, Writing - Original draft preparation, Writing - review $\&$ editing}

\affiliation[1]{organization={Shen Yuan Honors College, Beihang University},
                addressline={Xueyuan Road No.37}, 
                city={Beijing},
                postcode={100191}, 
                country={People's Republic of China}}

\affiliation[2]{organization={School of Astronautics, Beihang University},
                addressline={Xueyuan Road No.37}, 
                postcode={100191}, 
                city={Beijing},
                country={People's Republic of China}}

\author[2]{Shengping Gong}[style=chinese]
\cormark[1]
\ead{gongsp@buaa.edu.cn}
\credit{Conceptualization, Funding acquisition, Writing - review $\&$ editing}

\author[2]{Peng Shi}[style=chinese]
\ead{shipeng@buaa.edu.cn}
\credit{Methodology, Formal analysis, Writing - review $\&$ editing}

\cortext[cor1]{Corresponding author}

\begin{abstract}
This paper develops an analytical strategy for solving the linear quadratic pursuit-evasion game in arbitrary Keplerian reference orbits. The motion of the pursuer and evader is described using the controlled Tschauner-Hempel equations, and the optimal game strategies of the pursuer and evader are presented by the solution of the differential Riccati equation.The analytical solution of the differential Riccati equation is presented for elliptic, parabolic, and hyperbolic reference orbits, thereby enabling an analytical pursuit-evasion game strategy. Then, the procedure to solve the pursuit-evasion game using this analytical strategy is proposed. Simulations of pursuit-evasion game in elliptic, parabolic, and hyperbolic reference orbits validate the effectiveness of the developed analytical strategy. Results indicates that the analytical strategy saves the CPU time by more than 99.8$\%$ compared to the numerical one, highlighting the efficiency of the developed strategy. The developed analytical strategy is also applicable to pursuit-evasion game scenarios considering orbital disturbances. Compared to the conventional strategy, which succeed in only two out of six test scenarios, the developed strategy achieves success in all six cases, particularly demonstrating its effectiveness in high-eccentricity cases.
\end{abstract}

\begin{highlights}
\item The pursuit-evasion scenario is extended to arbitrary reference orbits.
\item Analytical strategy for pursuit-evasion games is derived.
\item Effectiveness of analytical strategy in pursuit-evasion games with disturbances are discussed.
\end{highlights}


\maketitle

\section{Introduction}\label{sec1}
Orbital pursuit-evasion game techniques play an important role in spacecraft interception and rendezvous missions \cite{prince2019elliptical}, as well as in the defense against space debris \cite{zhengqing2024space} and near-Earth asteroids (NEAs) \cite{vernis2018genetic}. Several scenarios of orbital pursuit-evasion games have been considered, such as two-player games \cite{pang2024solving} and multiple-player games \cite{li2024nash}. The two-player pursuit-evasion game, first introduced by Isaacs \cite{isaacs1999differential}, has attracted considerable attention from scholars. In this scenario, one player acts as the pursuer, aiming to intercept or rendezvous with the other player, the evader, who strives to avoid capture or interception \cite{jagat2017nonlinear}.

Solving the orbital pursuit-evasion game presents challenges due to strong orbital dynamics constraints, strong nonlinearities, and the high dimensionality of the problem \cite{pang2024solving}. To address these challenges, Pontani and Conway \cite{pontani2009numerical} developed a three-dimensional game and derived the necessary conditions for saddle-point solutions. Their works focused on the absolute motion of the pursuer and evader around the Earth. Meanwhile, orbital pursuit-evasion games based on reference orbits and relative dynamics have been further developed. Stupik et al. \cite{stupik2012optimal} explored optimal pursuit-evasion games in the Hill reference frame (i.e., the local-vertical/local-horizontal (LVLH) frame \cite{alfriend2009spacecraft}), introducing the Clohessy-Wiltshire (CW) equations \cite{clohessy1960terminal} to describe orbital pursuit-evasion games near a circular orbit. The CW equations have since facilitated rapid advancements in orbital pursuit-evasion game techniques. Tartaglia and Innocenti \cite{tartaglia2016game} investigated the linear quadratic (LQ) pursuit-evasion game using the CW relative dynamics. Li et al. \cite{li2019dimension} developed a dimension-reduction method for pursuit-evasion games in the context of the CW equations and further explored the saddle-point solutions in the J2-perturbed CW relative dynamics \cite{li2020saddle}. Ye et al. \cite{dong2018satellite,ye2020satellite} employed the CW equations to explore the pursuit-evasion game with vector guidance and different thrust configurations. The use of circular reference orbits and the CW equations transforms nonlinear pursuit-evasion games into linear ones, simplifying the solution of the problem. However, the CW equations become less accurate for describing the relative motion of the pursuer and evader in orbits with high eccentricity. To address this limitation, elliptic reference orbits and the corresponding Tschauner-Hempel (TH) equations \cite{dang2017solutions} have been introduced. Jagat and Sinclair \cite{jagat2017nonlinear} discussed the LQ pursuit-evasion game in nonlinear relative dynamics and the TH equations. Prince et al. \cite{prince2019elliptical} summarized proximity differential games in elliptic reference orbits, and Pang et al. \cite{pang2024solving} provided precise gradients used in optimization methods to rapidly solve pursuit-evasion games in elliptic reference orbits. Though previous works primarily focused on the pursuit-evasion game in circular or elliptic reference orbits, there is growing interest in extending these scenarios to parabolic and hyperbolic reference orbits. From a mathematical perspective, Keplerian orbits can be classified into elliptic (including circular), parabolic, and hyperbolic orbits \cite{battin1999an}. Investigating pursuit-evasion games in parabolic and hyperbolic reference orbits not only enhances the theoretical framework of orbital game theory but also offers broader and more universal scenarios for future research. Moreover, missions involving high-speed spacecraft, asteroid, and active debris defense are emerging as key areas of focus \cite{vernis2018genetic,zhengqing2024space}. These scenarios often involve high-speed targets, where relative dynamics in circular or elliptic reference orbits may no longer provide sufficient accuracy. Studies on pursuit-evasion games in parabolic and hyperbolic reference orbits could therefore serve as valuable references for proximity missions involving such high-speed objects. Moreover, pursuit-evasion games in parabolic and hyperbolic reference orbits are also characterized by unique challenges, such as increased fuel consumption and non-periodic orbital dynamics. These factors result in significant differences in solutions of pursuit-evasion games compared to those in circular and elliptic reference orbits. However, current research has yet to fully address these issues. To bridge this gap, this paper extends pursuit-evasion game scenarios to parabolic and hyperbolic reference orbits and develops a universal game strategy applicable to arbitrary Keplerian reference orbits.

For a orbital pursuit-evasion game, several methods have been proposed to solve the problem, including the LQ method \cite{li2024nash,li2024orbital}, methods based on Pontryagin Maximum Principle \cite{sandberg2022autonomous,raigoza2022autonomous,pang2024solving}, and artificial intelligence methods \cite{zhao2023prd,wilt2022microsatellite}. Among these, the LQ method, characterized by a quadratic cost function based on a linear model \cite{gong2008relative}, can provide a feedback control strategy and has been widely used in the aerospace fields \cite{anderson2007optimal}. Therefore, this paper focuses on the LQ pursuit-evasion game. In the LQ pursuit-evasion game, determining the control input reduces to solving the Riccati equation \cite{jagat2017nonlinear}. The Riccati equation can be categorized into the differential Riccati equation (DRE) and the algebraic Riccati equation (ARE) \cite{gong2008relative}, with the DRE being applicable to finite-time domain problems and the ARE to infinite-time domain problems. As this paper focuses on finite-time domain pursuit-evasion games, the DRE is adopted. However, since the DRE is a nonlinear equation, current works \cite{li2024nash,li2024orbital} investigating the LQ pursuit-evasion game mostly used the numerical or discretization methods to solve it. These methods can be further improved in terms of computational efficiency by analytical or semi-analytical solutions. Liao \cite{liao2021research} developed a semi-analytical solution of the DRE using the inverse Riccati differential equation for LQ pursuit-evasion games in circular reference orbits. However, this method requires the weighting matrix in the terminal condition to be strictly positive definite to ensure the invertibility of the Riccati matrix. Additionally, the solution is derived under the assumption of time-independence, which limits its applicability to time-dependent systems, such as the TH equations in elliptic, parabolic, and hyperbolic reference orbits. To address these limitations, using the state transition matrix (STM) developed by Dang \cite{dang2017new}, this paper derives an analytical solution of the DRE for the LQ pursuit-evasion game in arbitrary Keplerian reference orbits. Differing from Liao’s semi-analytical solution \cite{liao2021research} restricted to circular reference orbits, the analytical solution of the DRE is derived using the variation of parameters (VOP) method \cite{alfriend2009spacecraft} and variable transformation \cite{li2012fuel}, which are also applicable to the time-dependent TH equations. The obtained solution does not require the weighting matrix in the terminal condition to be strictly positive definite and is applicable to more scenarios than Liao's solution. Based on the obtained solution, an analytical game strategy is generated. A procedure using the analytical strategy is then proposed for solving pursuit-evasion games, and simulations demonstrate the effectiveness and efficiency of the developed analytical strategy. Furthermore, the effectiveness of the developed strategy in the cases with orbital disturbances is also analyzed, and comparison between the developed analytical and conventional strategies highlights the advantage of the developed strategy. Above all, the contributions of this paper can be summarized as follows:
\begin{enumerate}[label=(\arabic*)]
\item We extend the pursuit-evasion game scenarios to arbitrary Keplerian reference orbits, enhancing the mathematical framework of orbital game theory. The strategy and procedure developed in this paper can be extended to various scenarios, such as pursuing high-speed spacecraft or defending against high-speed space debris and NEAs.
\item An analytical solution of the DRE for the LQ pursuit-evasion game in arbitrary Keplerian reference orbits is derived based on the STM developed by Dang \cite{dang2017solutions} and the VOP method \cite{alfriend2009spacecraft}. The analytical game strategy is then presented using this solution. Comparisons between the results obtained from the analytical and numerical strategies demonstrate the effectiveness and efficiency of the analytical strategy.
\end{enumerate}

The rest of this paper is organized as follows. Section \ref{sec2} presents the problem statement, including the controlled TH equations and the definition of LQ pursuit-evasion game. Section \ref{sec3} derives the analytical game strategy in arbitrary Keplerian reference orbits and proposes the procedure using the analytical strategy. Section \ref{sec4} presents and discusses the results of the analytical strategy, highlighting its performance and effectiveness in various scenarios. Finally, conclusions are drawn in Section \ref{sec5}.

\section{Problem Statement}\label{sec2}
In this section, the pursuit-evasion game considered in this paper is stated, including the controlled TH equations and the definition of the LQ pursuit-evasion game.

\subsection{Controlled TH Equations}\label{subsec2.1}

Several linear models have been proposed to describe the relative motion, including the Clohessy-Wiltshire (CW) equations \cite{clohessy1960terminal} and Tschauner-Hempel (TH) equations \cite{dang2017solutions}. The CW equations are derived under the assumption of a circular reference orbit, while the TH equations are derived for arbitrary Keplerian reference orbits. Since the pursuit-evasion game is considered in arbitrary Keplerian reference orbits, the TH equations are used to describe the motion of the purser and evader. The equations are expressed in the local-vertical/local-horizontal (LVLH) frame \cite{alfriend2009spacecraft}. The relative state is defined as $\bm{X}=\left[x,\text{ }y,\text{ }z,\text{ }\dot{x},\text{ }\dot{y},\text{ }\dot{z}\right]^{\text{T}}$, where ${\dot{ \left(  \cdot  \right)} }$ denotes the time derivative. The controlled TH equations in the LVLH frame in terms of time derivatives can be expressed as \cite{prince2019elliptical}:

\begin{equation}
\left\{ \begin{gathered}
  \ddot x = 2\dot f\dot y + \ddot fy + {\left( {\dot f} \right)^2}x + 2\mu {\left( {\frac{{\dot f}}{h}} \right)^{3/2}}x + {u_x} \hfill \\
  \ddot y =  - 2\dot f\dot x - \ddot fx + {\left( {\dot f} \right)^2}y - \mu {\left( {\frac{{\dot f}}{h}} \right)^{3/2}}y + {u_y} \hfill \\
  \ddot z =  - \mu {\left( {\frac{{\dot f}}{h}} \right)^{3/2}}z + {u_z} \hfill \\ 
\end{gathered}  \right.\label{eq100}
\end{equation}
where $h=\sqrt{\mu{p}}$ and $\mu$ denotes the Earth gravitational constant. The parameters $p$ and $f$ denote the semilatus rectum and true anomaly of the reference orbit, respectively. To simplify the expression of the TH equations, the state is transformed into the following form \cite{yamanaka2002new}:

\begin{equation}
\left[ {\begin{array}{*{20}{c}}
  {\tilde x} \\ 
  {\tilde y} \\ 
  {\tilde z} 
\end{array}} \right] = \rho(f) \left[ {\begin{array}{*{20}{c}}
  x \\ 
  y \\ 
  z 
\end{array}} \right]\label{eq1}
\end{equation}
where $\rho(f)=1+e\cos{f}$. The parameter $e$ denotes the eccentricity of the reference orbit. Let ${\left(  \cdot  \right)' }$ denote the derivative with respect to true anomaly, the relationship between ${\left(  \cdot  \right)' }$ and ${\dot{ \left(  \cdot  \right)} }$ can be expressed as ($w$ is an arbitrary variable):

\begin{equation}
\begin{gathered}
  \dot w = \dot fw' = n{\rho ^2}w' \hfill \\
  \ddot w = \ddot fw' + {\left( {\dot f} \right)^2}w'' = {n^2}{\rho ^3}\left( {\rho w'' - 2e\sin fw'} \right) \hfill \\ 
\end{gathered} \label{eq101}
\end{equation}
where $n=\sqrt{\mu/p^3}$. Based on the relationship described in Eq. \eqref{eq101}, the controlled TH equations in terms of derivatives with respect to true anomaly are expressed as \cite{li2012fuel}:

\begin{equation}
\left\{ \begin{gathered}
  \tilde x'' - \frac{3}{\rho }\tilde x - 2\tilde y' = \frac{\beta }{{{\rho ^3}}}{u_x} \hfill \\
  \tilde y'' + 2\tilde x' = \frac{\beta }{{{\rho ^3}}}{u_y} \hfill \\
  \tilde z'' + \tilde z = \frac{\beta }{{{\rho ^3}}}{u_z} \hfill \\ 
\end{gathered}  \right.\label{eq2}
\end{equation}
where
\begin{equation}
\beta  = \frac{1}{{{n^2}}} = \frac{1}{{{\mu  \mathord{\left/
 {\vphantom {\mu  {{p^3}}}} \right.
 \kern-\nulldelimiterspace} {{p^3}}}}}\label{eq40}
\end{equation}
Denoting $\bm{u}=\left[u_x,\text{ }u_y,\text{ }u_z\right]^{\text{T}}$, Eq. \eqref{eq2} can be expressed in the following matrix form:

\begin{equation}
\tilde{\bm X}' = \bm{A}\tilde{\bm X} + \bm{Bu}\label{eq4}
\end{equation}
with

\begin{equation}
\bm{A} = \left[ {
  {\begin{array}{*{20}{c}}
  0&0&0&1&0&0 \\ 
  0&0&0&0&1&0 \\ 
  0&0&0&0&0&1 \\
  {\frac{3}{\rho }}&0&0&0&2&0 \\
  0&0&0&{ - 2}&0&0 \\
  0&0&{ - 1}&0&0&0
\end{array}}} \right]\label{eq5}
\end{equation}

\begin{equation}
\bm{B} = \frac{\beta }{{{\rho ^3}}}\left[ {\begin{array}{*{20}{c}}
  {{\bm{O}_{3 \times 3}}} \\ 
  {{\bm{I}_{3 \times 3}}} 
\end{array}} \right]\label{eq6}
\end{equation}
When the $e$ is set to 0 in the equations, the TH equations is simplified into the CW equations \cite{yamanaka2002new}. Based on the controlled TH equations, the LQ pursuit-evasion game is defined.

\subsection{LQ Pursuit-Evasion Game}\label{subsec2.2}

To derive the analytical strategy of the pursuit-evasion game, we focus on a simple scenario, i.e., a two-player game. In two-player pursuit-evasion games, the problem can be categorized into impulsive pursuit-evasion games \cite{ma2024delta,xie2024game} and continuous thrust pursuit-evasion games \cite{prince2019elliptical,pang2024solving}. In this paper, we consider a scenario of the continuous thrust pursuit-evasion game, specifically the LQ pursuit-evasion game, and investigate its analytical strategy. The motion of the purser and evader is described by the controlled TH equations presented in Eq. \eqref{eq4}:

\begin{equation}
\left\{ \begin{gathered}
  {{\tilde{\bm X}'}_p} = \bm{A}{{\tilde{\bm X}}_p} + \bm{B}{\bm{u}_p} \hfill \\
  {{\tilde{\bm X}'}_e} = \bm{A}{{\tilde{\bm X}}_e} + \bm{B}{\bm{u}_e} \hfill \\ 
\end{gathered}  \right.\label{eq7}
\end{equation}
where the subscript ‘\textit{p}’ denotes the quantities associated with the purser, while the subscript ‘\textit{e}’ denotes the quantities associated with the evader. Therefore, defining $\tilde{\bm x} = {\tilde{\bm X}_p} - {\tilde{\bm X}_e}$, the dynamical equations between the purser and evader can be expressed as:

\begin{equation}
\tilde{\bm x}' = \bm{A}\tilde{\bm x} + \bm{B}{\bm{u}_p} - \bm{B}{\bm{u}_e}\label{eq8}
\end{equation}

In the LQ pursuit-evasion game, the cost function to be minimized is expressed as \cite{liao2021research}:

\begin{equation}
J = \frac{1}{2}\left( {\tilde{\bm x}{{\left( {{f_f}} \right)}^{\text{T}}}\bm{S}\tilde{\bm x}\left( {{f_f}} \right)} \right) + \frac{1}{2}\int_{f_0}^{{f_f}} {\left( {{\bm{u}_p}^{\text{T}}{\bm{R}_p}{\bm{u}_p} - {\bm{u}_e}^{\text{T}}{\bm{R}_e}{\bm{u}_e}} \right){\text{d}}f} \label{eq9}
\end{equation}
where $\bm{S}$, $\bm{R}_p$, and $\bm{R}_e$ denote the weighting matrices. In this paper, the intercept game is considered, and the weighting matrices are set as $\bm{S} = {\text{diag}}\left( s_r{{\bm{I}_{3 \times 3}},{\text{ }}{\bm{O}_{3 \times 3}}} \right)$, ${\bm{R}_p} = {r_p}{\bm{I}_{3 \times 3}}$, and ${\bm{R}_e} = {r_e}{\bm{I}_{3 \times 3}}$. The $f_0$ and $f_f$ indicate the initial and terminal true anomalies of the pursuit-evasion game. Therefore, the Hamiltonian function is defined as \cite{liao2021research}:

\begin{equation}
\mathcal{H} = \frac{1}{2}\left( {{\bm{u}_p}^{\text{T}}{\bm{R}_p}{\bm{u}_p} - {\bm{u}_e}^{\text{T}}{\bm{R}_e}{\bm{u}_e}} \right) + {\bm{\lambda }^{\text{T}}}\left( {\bm{A}\tilde{\bm x} + \bm{B}{\bm{u}_p} - \bm{B}{\bm{u}_e}} \right)\label{eq10}
\end{equation}
where $\bm{\lambda }={\left[ {{\lambda _{\tilde x}},{\text{ }}{\lambda _{\tilde y}},{\text{ }}{\lambda _{\tilde z}},{\text{ }}{\lambda _{\tilde x'}},{\text{ }}{\lambda _{\tilde y'}},{\text{ }}{\lambda _{\tilde z'}}} \right]^{\text{T}}}$ denotes the costate associated with $\tilde{\bm x}$. The terminal functional is defined as \cite{liao2021research}:
\begin{equation}
\Phi  = \frac{1}{2}{\tilde{\bm x}^{\text{T}}}\left( {{f_f}} \right)\bm{S}\tilde{\bm x}\left( {{f_f}} \right)\label{eq18}
\end{equation}
Therefore, the equations describing the pursuit-evasion game are presented as follows:
\begin{equation}
\left\{\begin{gathered}
  \tilde{\bm x}' = \frac{{\partial \mathcal{H}}}{{\partial \bm{\lambda }}} = \bm{A}\tilde{\bm x} + \bm{B}{\bm{u}_p} - \bm{B}{\bm{u}_e} \hfill \\
  \bm{\lambda }' =  - \frac{{\partial \mathcal{H}}}{{\partial \tilde{\bm x}}} =  - {\bm{A}^{\text{T}}}\bm{\lambda } \hfill \\ 
\end{gathered} \right. \label{eq11}
\end{equation}
The terminal condition is expressed as:
\begin{equation}
\bm{\lambda }\left( {{f_f}} \right) = \frac{{\partial \Phi }}{{\partial \tilde{\bm x}}}|_{\tilde{\bm x}=\tilde{\bm x}\left(f_f\right)} = \bm{S}\tilde{\bm x}\left( {{f_f}} \right)\label{eq19}
\end{equation}
According to the optimal necessary conditions of the differential game \cite{engwerda2005lq}, the optimal control input of the LQ pursuit-evasion game can be calculated by:
\begin{equation}
\left\{ \begin{gathered}
  {\bm{u}_p} =  - {\left( {{\bm{R}_p}} \right)^{ - 1}}{\bm{B}^{\text{T}}}\bm{\lambda } \hfill \\
  {\bm{u}_e} =  - {\left( {{\bm{R}_e}} \right)^{ - 1}}{\bm{B}^{\text{T}}}\bm{\lambda } \hfill \\ 
\end{gathered}  \right.\label{eq12}
\end{equation}

Substituting Eq. \eqref{eq12} for Eq. \eqref{eq11}, we can obtain the following equations with the terminal condition Eq. \eqref{eq19}:

\begin{equation}
\tilde{\bm x}' = \bm{A}\tilde{\bm x} - \bm{B}{\left( {{\bm{R}_p}} \right)^{ - 1}}{\bm{B}^{\text{T}}}\bm{\lambda } + \bm{B}{\left( {{\bm{R}_e}} \right)^{ - 1}}{\bm{B}^{\text{T}}}\bm{\lambda }\label{eq13}
\end{equation}
Therefore, Eq. \eqref{eq11} can be restructured in the following matrix form:
\begin{equation}
\left[ {\begin{array}{*{20}{c}}
  {\tilde{\bm x}'} \\ 
  {\bm{\lambda }'} 
\end{array}} \right] = \left[ {\begin{array}{*{20}{c}}
  \bm{A}&{ - \bm{B}\left( {{{\left( {{\bm{R}_p}} \right)}^{ - 1}} - {{\left( {{\bm{R}_e}} \right)}^{ - 1}}} \right){\bm{B}^{\text{T}}}} \\ 
  \bm{O}&{ - {\bm{A}^{\text{T}}}} 
\end{array}} \right]\left[ {\begin{array}{*{20}{c}}
  {\bm{\tilde x}} \\ 
  \bm{\lambda } 
\end{array}} \right]\label{eq14}
\end{equation}
Then, the solution of Eq. \eqref{eq14} is expressed as:
\begin{equation}
\left[ {\begin{array}{*{20}{c}}
  {\tilde{\bm x}\left( {{f}} \right)} \\ 
  {\bm{\lambda }\left( {{f}} \right)} 
\end{array}} \right] = \left[ {\begin{array}{*{20}{c}}
  {{\bm{\Omega }_{11}}\left( {{f},{\text{ }}{f_0}} \right)}&{{\bm{\Omega }_{12}}\left( {{f},{\text{ }}{f_0}} \right)} \\ 
  {{\bm{\Omega }_{21}}\left( {{f},{\text{ }}{f_0}} \right)}&{{\bm{\Omega }_{22}}\left( {{f},{\text{ }}{f_0}} \right)} 
\end{array}} \right]\left[ {\begin{array}{*{20}{c}}
  {\tilde{\bm x}\left( {{f_0}} \right)} \\ 
  {\bm{\lambda }\left( {{f_0}} \right)} 
\end{array}} \right]\label{eq15}
\end{equation}
where ${{\bm{\Omega }_{11}}\left( {{f},{\text{ }}{f_0}} \right)}=\exp\left(\int_{{f_0}}^{{f}} \bm{A} {\text{ d}}f\right)$ and ${{\bm{\Omega }_{22}}\left( {{f},{\text{ }}{f_0}} \right)}={{\bm{\Omega }_{11}}^{-\text{T}}\left( {{f},{\text{ }}{f_0}} \right)}$. ${{\bm{\Omega }_{21}}\left( {{f},{\text{ }}{f_0}} \right)}$ is a zero matrix, and ${{\bm{\Omega }_{12}}\left( {{f},{\text{ }}{f_0}} \right)}$ can be calculated by the variation
of parameters (VOP) method \cite{li2012fuel}:
\begin{align}
  \tilde{\bm x}\left( {{f}} \right) & = {\bm{\Omega }_{11}}\left( {{f},{\text{ }}{f_0}} \right)\tilde{\bm x}\left( {{f_0}} \right) + \int_{{f_0}}^{{f}} {{\bm{\Omega }_{11}}\left( {{f},{\text{ }}\theta} \right)} \left( { - \bm{B}\left( {{{\left( {{\bm{R}_p}} \right)}^{ - 1}} - {{\left( {{\bm{R}_e}} \right)}^{ - 1}}} \right){\bm{B}^{\text{T}}}} \right)\bm{\lambda }\left( \theta \right){\text{d}}\theta \\ 
  \notag & = {\bm{\Omega }_{11}}\left( {{f},{\text{ }}{f_0}} \right)\tilde{\bm x}\left( {{f_0}} \right) + \int_{{f_0}}^{{f}} {{\bm{\Omega }_{11}}\left( {{f},{\text{ }}\theta} \right)} \left( { - \bm{B}\left( {{{\left( {{\bm{R}_p}} \right)}^{ - 1}} - {{\left( {{\bm{R}_e}} \right)}^{ - 1}}} \right){\bm{B}^{\text{T}}}} \right){\bm{\Omega }_{22}}\left( {\theta,{\text{ }}{f_0}} \right){\text{d}}\theta \cdot \bm{\lambda }\left( {{f_0}} \right) \\ 
  \notag & = {\bm{\Omega }_{11}}\left( {{f},{\text{ }}{f_0}} \right)\tilde{\bm x}\left( {{f_0}} \right) + \int_{{f_0}}^{{f}} {{\bm{\Omega }_{11}}\left( {{f},{\text{ }}\theta} \right)} \left( { - \bm{B}\left( {{{\left( {{\bm{R}_p}} \right)}^{ - 1}} - {{\left( {{\bm{R}_e}} \right)}^{ - 1}}} \right){\bm{B}^{\text{T}}}} \right){\bm{\Omega }_{11}}^{ - {\text{T}}}\left( {\theta,{\text{ }}{f_0}} \right){\text{d}}\theta \cdot \bm{\lambda }\left( {{f_0}} \right) \label{eq16} 
\end{align}

\begin{equation}
{\bm{\Omega }_{12}}\left( {{f},{\text{ }}{f_0}} \right) = \int_{{f_0}}^{{f}} {{\bm{\Omega }_{11}}\left( {{f},{\text{ }}\theta} \right)} \left( { - \bm{B}\left( {{{\left( {{\bm{R}_p}} \right)}^{ - 1}} - {{\left( {{\bm{R}_e}} \right)}^{ - 1}}} \right){\bm{B}^{\text{T}}}} \right){\bm{\Omega }_{11}}^{ - {\text{T}}}\left( {\theta,{\text{ }}{f_0}} \right){\text{d}}\theta\label{eq17}
\end{equation}
According to Eq. \eqref{eq19} and Eqs. \eqref{eq15}-\eqref{eq17}, a linear relationship between $\bm{\tilde x}$ and $\bm{\lambda}$ can be observed. It is assumed that:
\begin{equation}
\bm{\lambda }\left( f \right) = \bm{P}\left( f \right)\tilde{\bm x}\left( f \right)\label{eq20}
\end{equation}
Substituting Eq. \eqref{eq20} into Eq. \eqref{eq14}, we can obtain the following equation:
\begin{equation}
\left[ {\bm{P}'\left( f \right) + {\bm{A}^{\text{T}}}\bm{P}\left( f \right) + \bm{P}\left( f \right)\bm{A} - \bm{P}\left( f \right)\bm{B}\left( {{{\left( {{\bm{R}_p}} \right)}^{ - 1}} - {{\left( {{\bm{R}_e}} \right)}^{ - 1}}} \right){\bm{B}^{\text{T}}}\bm{P}\left( f \right)} \right]\tilde{\bm x}\left( f \right) = \textbf{0}\label{eq21}
\end{equation}
Equation \eqref{eq21} holds for any state during the pursuit-evasion game. Therefore, Eq. \eqref{eq21} is transformed into the following DRE \cite{liao2021research}:

\begin{equation}
\bm{P}'\left( f \right)\bm{ = } - {\bm{A}^{\text{T}}}\bm{P}\left( f \right) - \bm{P}\left( f \right)\bm{A} + \bm{P}\left( f \right)\bm{B}\left( {{{\left( {{\bm{R}_p}} \right)}^{ - 1}} - {{\left( {{\bm{R}_e}} \right)}^{ - 1}}} \right){\bm{B}^{\text{T}}}\bm{P}\left( f \right)\label{eq22}
\end{equation}

\begin{equation}
\bm{P}\left( {{f_f}} \right) = \bm{S}\label{eq23}
\end{equation}

In this paper, the DRE is adopted to generate the optimal control input of the pursuer and evader. According to Eq. \eqref{eq12} and Eq. \eqref{eq20}, the control input generated from the solution of the DRE Eq. \eqref{eq22} is a feedback control. Generally, the DRE requires numerical backward-in-time integration \cite{li2024nash}. However, based on the universal solutions for the TH equations developed by Ref. \cite{dang2017solutions}, an analytical solution of the DRE can be obtained. Consequently, the analytical pursuit-evasion game strategy is developed for arbitrary Keplerian reference orbits.

\section{Analytical Pursuit-Evasion Game Strategy}\label{sec3}
In this section, the analytical pursuit-evasion game strategy is derived for elliptic, parabolic, and hyperbolic reference orbits based on the DRE. Subsequently, a procedure using the derived analytical strategy for solving the pursuit-evasion game is proposed.
\subsection{Solution of DRE}\label{subsec3.1}
According to Eq. \eqref{eq15} and Eq. \eqref{eq20}, when the matrix $\left( {{\bm{\Omega }_{22}}\left( {{f_f},{f}} \right) - \bm{S}{\bm{\Omega }_{12}}\left( {{f_f},{f}} \right)} \right)$ is invertible, the analytical solution of the DRE can be expressed as:

\begin{equation}
  \bm{P}{\left(f\right)} = {\left( {{\bm{\Omega }_{22}}\left( {{f_f},{f}} \right) - \bm{S}{\bm{\Omega }_{12}}\left( {{f_f},{f}} \right)} \right)^{ - 1}} \cdot \bm{S}{\bm{\Omega }_{11}}\left( {{f_f},{f}} \right)  
 \label{eq1000}
\end{equation}
According to Eq. \eqref{eq1000}, the analytical solution of the DRE Eq. \eqref{eq22} involves ${\bm{\Omega }_{11}}\left( {{f_f},{f}}\right)$, ${\bm{\Omega }_{12}}\left( {{f_f},{f}}\right)$, and ${\bm{\Omega }_{22}}\left( {{f_f},{f}}\right)$. Therefore, according to Eq. \eqref{eq17}, the state transition matrix (STM) for relative motion (i.e., ${{\bm{\Omega }_{11}}\left( {{f_2},{\text{ }}{f_1}} \right)}$ from $f_1$ to $f_2$ in the reference orbit) \cite{dang2017solutions} is an effective tool for obtaining the analytical solution of the DRE Eq. \eqref{eq22}. Since this paper focuses on the pursuit-evasion game in arbitrary Keplerian reference orbits, the STM for relative motion in arbitrary Keplerian reference orbits is required. Carter \cite{carter1998state} first developed a solution of the TH equations in arbitrary Keplerian reference orbits and summarized it in the form of an STM. Yamanaka and Anderson \cite{yamanaka2002new} derived a new closed-form STM of the elliptic TH equations based on Carter's work, enhancing engineering applications \cite{li2012fuel,pang2024solving}. Extending the works of Carter, Yamanaka, and Anderson, Dang \cite{dang2017new} developed a new solution (i.e., a new STM) of the TH equations and proved the equivalence between the previous two STMs \cite{carter1998state,yamanaka2002new}. This paper employs the new STM developed by Dang \cite{dang2017solutions} for its elegant form and derives the analytical solution of the DRE in arbitrary Keplerian reference orbits. According to Ref. \cite{dang2017solutions}, ${{\bm{\Omega }_{11}}\left( {{f_2},{\text{ }}{f_1}} \right)}$ is expressed as (in this Section, $f_1$ and $f_2$ are used to denote any initial and terminal true anomaly of the STMs, e.g., $f_1$ and $f_2$ can represent the present true anomaly $f$ and terminal true anomaly $f_f$ in the pursuit-evasion game, respectively):

\begin{equation}
{\bm{\Omega }_{11}}\left( {{f_2},{\text{ }}{f_1}} \right) = \phi \left( {{f_2}} \right){\phi ^{ - 1}}\left( {{f_1}} \right)\label{eq24}
\end{equation}
where
\begin{equation}
\phi \left( f \right) = \left[ {\begin{array}{*{20}{c}}
  {\begin{array}{*{20}{c}}
  {{\varphi _1}}&{{\varphi _2}}&{{\varphi _3}}&0&0&0 \\ 
  { - 2S\left( {{\varphi _1}} \right)}&{ - 2S\left( {{\varphi _2}} \right)}&{ - S\left( {2{\varphi _3} + 1} \right)}&1&0&0 \\ 
  0&0&0& 0&{\cos f}&{\sin f} \\
  {{\varphi _1}' }&{{\varphi _2}' }&{{\varphi _3}' }&0&0&0 \\
   { - 2{\varphi _1}}&{ - 2{\varphi _2}}&{ - 2{\varphi _3} - 1}&0&0&0 \\
   0&0&0&0&{ - \sin f}&{\cos f} 
\end{array}} 
\end{array}} \right]\label{eq28}
\end{equation}

\begin{equation}
{\phi ^{ - 1}}\left( f \right) = \left[ {\begin{array}{*{20}{c}}
  {\begin{array}{*{20}{c}}
  {4S\left( {{\varphi _2}} \right) + {\varphi _2}' }&0&0&{ - {\varphi _2}}&{2S\left( {{\varphi _2}} \right)}&0 \\ 
  { - 4S\left( {{\varphi _1}} \right) - {\varphi _1}' }&0&0&{{\varphi _1}}&{ - 2S\left( {{\varphi _1}} \right)}&0 \\ 
  { - 2}&0&0 &0&{ - 1}&0 \\
   { - 2S\left( {2{\varphi _3} + 1} \right) - {\varphi _3}' }&1&0&{{\varphi _3}}&{ - S\left( {2{\varphi _3} + 1} \right)}&0 \\
   0&0&{\cos f}&0&0&{ - \sin f} \\
   0&0&{\sin f}&0&0&{\cos f} 
\end{array}} 
\end{array}} \right]\label{eq29}
\end{equation}

\begin{equation}
\left\{ \begin{gathered}
  {\varphi _1} = \rho \left( f \right)\sin f \hfill \\
  {\varphi _1}'  = \rho \left( f \right)\cos f - e{\sin ^2}f \hfill \\
  S\left( {{\varphi _1}} \right) =  - \cos f - \frac{1}{2}e{\cos ^2}f \hfill \\ 
\end{gathered}  \right.\label{eq30}
\end{equation}
For elliptic and hyperbolic reference orbits $\left( {e \ne 1} \right)$:
\begin{equation}
\left\{ \begin{gathered}
  {\varphi _2} = \frac{{e{\varphi _1}}}{{1 - {e^2}}}\left( {D\left( f \right) - 3eL\left( f \right)} \right) - \frac{{\cos f}}{{\rho \left( f \right)}} \hfill \\
  {\varphi _3} =  - \frac{{{\varphi _1}}}{{1 - {e^2}}}\left( {D\left( f \right) - 3eL\left( f \right)} \right) - \frac{{{{\cos }^2}f}}{{\rho \left( f \right)}} - {\cos ^2}f \hfill \\
  {\varphi _2}'  = \frac{{e{\varphi _1}' }}{{1 - {e^2}}}\left( {D\left( f \right) - 3eL\left( f \right)} \right) + \frac{{e\sin f\cos f}}{{{\rho ^2}\left( f \right)}} + \frac{{\sin f}}{{\rho \left( f \right)}} \hfill \\
  {\varphi _3}'  = 2\left( {{\varphi _1}' S\left( {{\varphi _2}} \right) - {\varphi _2}' S\left( {{\varphi _1}} \right)} \right) \hfill \\
  S\left( {{\varphi _2}} \right) =  - \frac{{{\rho ^2}\left( f \right)\left( {D\left( f \right) - 3eL\left( f \right)} \right)}}{{2\left( {1 - {e^2}} \right)}} \hfill \\
  S\left( {2{\varphi _3} + 1} \right) = \frac{{e\sin f\left( {2 + e\cos f} \right)}}{{1 - {e^2}}} - \frac{{3{\rho ^2}\left( f \right)L\left( f \right)}}{{1 - {e^2}}} \hfill \\
  D\left( f \right) = \frac{{\sin f\left( {2 + e\cos f} \right)}}{{{\rho ^2}\left( f \right)}} \hfill \\
  L\left( f \right) = \int {\frac{1}{{{\rho ^2}\left( f \right)}}} {\text{d}}f = nt \hfill \\ 
\end{gathered}  \right.\label{eq31}
\end{equation}
For the parabolic reference orbit $\left( {e = 1} \right)$:
\begin{equation}
\left\{ \begin{gathered}
  {\varphi _2} = 2{\varphi _1}C\left( f \right) - \frac{{\cos f}}{{\rho \left( f \right)}} \hfill \\
  {\varphi _3} =  - 2{\varphi _1}C\left( f \right) - \frac{{{{\cos }^2}f}}{{\rho \left( f \right)}} - {\cos ^2}f \hfill \\
  {\varphi _2}'  = 2{\varphi _1}' C\left( f \right) + \frac{{\sin f\cos f}}{{{\rho ^2}\left( f \right)}} + \frac{{\sin f}}{{\rho \left( f \right)}} \hfill \\
  {\varphi _3}'  = 2\left( {{\varphi _1}' S\left( {{\varphi _2}} \right) - {\varphi _2}' S\left( {{\varphi _1}} \right)} \right) \hfill \\
  S\left( {{\varphi _2}} \right) =  - {\rho ^2}\left( f \right)C\left( f \right) \hfill \\
  S\left( {2{\varphi _3} + 1} \right) = 2\left( {{\rho ^2}\left( f \right)C\left( f \right) - \sin f} \right) - \sin f\cos f \hfill \\
  C\left( f \right) = \frac{1}{4}\tan \left( {\frac{f}{2}} \right) - \frac{1}{{20}}{\tan ^5}\left( {\frac{f}{2}} \right) \hfill \\ 
\end{gathered}  \right.\label{eq32}
\end{equation}
Therefore, ${{\bm{\Omega }_{22}}\left( {{f_2},{\text{ }}{f_1}} \right)}={{\bm{\Omega }_{11}}^{-\text{T}}\left( {{f_2},{\text{ }}{f_1}} \right)}$ can be expressed as:
\begin{equation}
{\bm{\Omega }_{22}}\left( {{f_2},{\text{ }}{f_1}} \right) = {\left( {{\phi ^{ - 1}}\left( {{f_2}} \right)} \right)^{\text{T}}}{\left( {\phi \left( {{f_1}} \right)} \right)^{\text{T}}}\label{eq25}
\end{equation}
Substituting Eqs. \eqref{eq24}-\eqref{eq25} into Eq. \eqref{eq17}, ${{\bm{\Omega }_{12}}\left( {{f_2},{\text{ }}{f_1}} \right)}$ can be calculated by:
\begin{equation}
{\bm{\Omega }_{12}}\left( {{f_2},{\text{ }}{f_1}} \right) = \frac{1}{{{n^4}}}\left( {\frac{1}{{{r_e}}} - \frac{1}{{{r_p}}}} \right)\bm{J}_1\label{eq26}
\end{equation}
where
\begin{align}\label{eq27}
\bm{J}_1 &= \phi \left( f_2 \right)\int_{{f_1}}^{f_2} {{\phi ^{ - 1}}\left( \theta  \right)} \left[ {\begin{array}{*{20}{c}}
  \bm{O}&\bm{O} \\ 
  \bm{O}&{\frac{1}{{{\rho ^6}(\theta)}}{\bm{I}_{3 \times 3}}} 
\end{array}} \right]{\phi ^{ - {\text{T}}}}\left( \theta  \right){\text{d}}\theta \cdot {\phi ^{\text{T}}}\left( {{f_1}} \right) \\
\notag & =\phi \left( f_2\right)\left({\bm{J}\left(f_2\right)-\bm{J}\left(f_1\right)}\right){\phi ^{\text{T}}}\left( {{f_1}} \right)
\end{align}
To evaluate the solution of the DRE in arbitrary Keplerian reference orbits, the specific expressions of the components in the matrix $\bm{J}$ are derived in three cases, i.e., elliptic $\left(0 \leqslant e < 1\right)$, parabolic $\left(e = 1\right)$, and hyperbolic $\left(e > 1\right)$ reference orbits based on the geometric relationships of the orbital elements \cite{curtis2020orbital}.

\subsubsection{Elliptic Reference Orbit (Case I)}\label{subsubsec3.1.1}
In the original expressions of the matrix $\bm{J}_1$ for an elliptic reference orbit, the integration cannot be performed directly due to the presence of $1/{\rho^i(f)}$. Therefore, the eccentric anomaly ($E$) is introduced to make Eq. \eqref{eq27} integrable \cite{li2012fuel}. The specific geometric relationships required are provided as follows \cite{curtis2020orbital}:
\begin{equation}
\left\{ \begin{gathered}
  \sin f = \frac{{\sqrt {1 - {e^2}} \sin E}}{{1 - e\cos E}} \hfill \\
  \cos f = \frac{{\cos E - e}}{{1 - e\cos E}} \hfill \\
  L\left( f \right) = nt = \frac{{E - e\sin E}}{{\sqrt {{{\left( {1 - {e^2}} \right)}^3}} }} \hfill \\
  {\text{d}}f = \frac{{\sqrt {1 - {e^2}} }}{{1 - e\cos E}}{\text{d}}E \hfill \\ 
\end{gathered}  \right.\label{eq33}
\end{equation}
The eccentric anomaly can be calculated by the following expressions, achieved by the Matlab® ‘atan2’ command:
\begin{equation}
E = {\text{atan2}}\left( {\frac{{\sqrt {1 - {e^2}} \sin f}}{{\rho \left( f \right)}},{\text{ }}\frac{{\cos f + e}}{{\rho \left( f \right)}}} \right) + f - {\text{atan2}}\left( {\sin f,{\text{ }}\cos f} \right)\label{eq34}
\end{equation}
Since ${\text{atan2}}\left( {\frac{{\sqrt {1 - {e^2}} \sin f}}{{\rho \left( f \right)}},{\text{ }}\frac{{\cos f + e}}{{\rho \left( f \right)}}} \right) \in \left(-\pi,\text{ }\pi\right]$ while $f\in\mathbb{R}$, to match the values of $E$ with the range of $f$ (i.e., $\left| {f - E} \right| < 2\pi $), the term $f - {\text{atan2}}\left( {\sin f,{\text{ }}\cos f} \right)$ is included when calculating the value of $E$. Then, Eq. \eqref{eq32} is substituted into Eq. \eqref{eq27}, and the integration is performed as follows:

\begin{align}\label{eqA}
  {\bm{J}_1} &= \phi \left( f_2 \right)\int_{{E_1}}^{E_2} {{{\hat \phi }^{ - 1}}\left( \vartheta  \right)} \left[ {\begin{array}{*{20}{c}}
  \bm{O}&\bm{O} \\ 
  \bm{O}&{\frac{1}{{{{\hat \rho }^6}\left( \vartheta  \right)}}{\bm{I}_{3 \times 3}}} 
\end{array}} \right]{{\hat \phi }^{ - {\text{T}}}}\left( \vartheta  \right)\frac{{\sqrt {1 - {e^2}} }}{{1 - e\cos \vartheta }}{\text{d}}\vartheta  \cdot {\phi ^{\text{T}}}\left( {{f_1}} \right) \\ 
   \notag &= \phi \left( f_2 \right)\left( {\bm{J}\left( f_2 \right) - \bm{J}\left( {{f_1}} \right)} \right){\phi ^{\text{T}}}\left( {{f_1}} \right) \\ 
   \notag & = \phi \left( f_2 \right)\left( {\hat{\bm J}\left( E_2 \right) - \hat{\bm J}\left( {{E_1}} \right)} \right){\phi ^{\text{T}}}\left( {{f_1}} \right)
\end{align}
where the superscript $\hat {\left(  \cdot  \right)}$ denotes the transformed function from $f$ to $E$. The specific expressions of the components in the matrix $\bm{J}$ can be found in Appendix \ref{secA1}. Notably, the expressions take the same form for the circular reference orbit ($e=0$). Since the derived analytical solution involves the terms $1/(1-e^2)^j$, which may lead to inaccuracies when solving the DRE for reference orbits with eccentricities close to 1, it is recommended to select reference orbits with $0 \leqslant e \leqslant 0.8$ when applying this analytical solution to elliptic reference orbits. In cases involving pursuit-evasion games in high-eccentricity elliptic reference orbits (i.e., $e>0.8$), it is recommended to transform the initial states of the pursuer and evader into those applicable to the TH equations in parabolic reference orbits and to use the analytical solution applicable to parabolic reference orbits, as detailed in Section \ref{subsubsec3.1.2} and Appendix \ref{secA2}. 
\subsubsection{Parabolic Reference Orbit (Case II)}\label{subsubsec3.1.2}
When $e=1$, the original expressions of the matrix $\bm{J}_1$ shown in Eq. \eqref{eq27} can be directly integrated with respect to the true anomaly $f$. The specific expressions of the components in the matrix $\bm{J}$ can be found in Appendix \ref{secA2}.
\subsubsection{Hyperbolic Reference Orbit (Case III)}\label{subsubsec3.1.3}
Similar to Case I, the hyperbolic eccentric anomaly $H$ is introduced to make Eq. \eqref{eq27} integrable. The required geometric relationships are provided as follows \cite{curtis2020orbital}:
\begin{equation}
\left\{ \begin{gathered}
  \sin f = \frac{{\sqrt {{e^2} - 1} \sinh H}}{{e\cosh H - 1}} \hfill \\
  \cos f = \frac{{e - \cosh H}}{{e\cosh H - 1}} \hfill \\
  L\left( f \right) = nt = \frac{{e\sinh H - H}}{{\sqrt {{{\left( {{e^2} - 1} \right)}^3}} }} \hfill \\
  {\text{d}}f = \frac{{\sqrt {{e^2} - 1} }}{{e\cosh H - 1}}{\text{d}}H \hfill \\ 
\end{gathered}  \right.\label{eq35}
\end{equation}
where
\begin{equation}
H = {\text{atanh}}\left( {\frac{{\sinh H}}{{\cosh H}}} \right) = {\text{atanh}}\left( {\frac{{\sqrt {{e^2} - 1} \sin f}}{{e + \cos f}}} \right)\label{eq36}
\end{equation}
Then, Eq. \eqref{eq32} is substituted into Eq. \eqref{eq27}, and the integration is performed as follows:

\begin{align}\label{eqB}
  {\bm{J}_1} & = \phi \left( f_2 \right)\int_{{H_1}}^{H_2} {{{\bar \phi }^{ - 1}}\left( \psi  \right)} \left[ {\begin{array}{*{20}{c}}
  \bm{O}&\bm{O} \\ 
  \bm{O}&{\frac{1}{{{{\bar \rho }^6}\left( \psi  \right)}}{\bm{I}_{3 \times 3}}} 
\end{array}} \right]{{\bar \phi }^{ - {\text{T}}}}\left( \psi  \right)\frac{{\sqrt {{e^2} - 1} }}{{e\cosh \psi  - 1}}{\text{d}}\psi  \cdot {\phi ^{\text{T}}}\left( {{f_1}} \right) \\ 
  \notag & = \phi \left( f_2 \right)\left( {\bm{J}\left( f_2 \right) - \bm{J}\left( {{f_1}} \right)} \right){\phi ^{\text{T}}}\left( {{f_1}} \right) \\ 
  \notag & = \phi \left( f_2 \right)\left( {\bar{\bm J}\left( H_2 \right) - \bar{\bm{J}}\left( {{H_1}} \right)} \right){\phi ^{\text{T}}}\left( {{f_1}} \right)  
\end{align}

where the superscript $\bar {\left( \cdot \right)}$ denotes the transformed function from $f$ to $H$. The specific expressions of the components in the matrix $\bm{J}$ can be found in Appendix \ref{secA3}. Similar to the case in the elliptic reference orbit, the terms $1/(e^2-1)^k$ may lead to inaccuracies when evaluating the analytical solution of the DRE in reference orbits with eccentricities close to 1. Therefore, it is recommended to select reference orbits with $e \geqslant 1.2$ when applying the analytical solution in hyperbolic reference orbits. When considering the pursuit-evasion game in reference orbits with $1 < e < 1.2$, it is recommended to transform the scenario into a parabolic reference orbit.

\subsubsection{Time-to-go Estimation
}\label{subsubsec3.1.4}

Since the boundary condition of the DRE presented in Eq. \eqref{eq23} is the value of matrix $P$ at the end of the pursuit-evasion game (i.e., at $f_f$), a time-to-go estimation (i.e., an estimation of the remaining true anomaly from the present true anomaly $f$ to the terminal true anomaly $f_f$ in the reference orbit) is required. This paper focuses on the intercept game. According to Ref. \cite{liao2021research}, the time-to-go is estimated as follows:

\begin{equation}
{f_f} - f \approx  -{\frac{{{{\tilde x}^2} + {{\tilde y}^2} + {{\tilde z}^2}}}{{\tilde x\tilde x' + \tilde y\tilde y' + \tilde z\tilde z'}}} \label{eq37}
\end{equation}
Then, the solution of the DRE Eq. \eqref{eq22} can be obtained from the Eq. \eqref{eq1000} and STMs mentioned in Sections \ref{subsubsec3.1.1}-\ref{subsubsec3.1.3} and Appendices \ref{secA1}-\ref{secA3}. Subsequently, the control inputs of the pursuer and evader in the LQ pursuit-evasion game, $\bm{u}_p$ and $\bm{u}_e$, can be calculated analytically using Eq. \eqref{eq12}. Therefore, in this paper, the analytical strategy is developed based on the analytical solution of the DRE \eqref{eq22}. 

\noindent \textbf{Remark 1}: Compared to the game strategy based on the semi-analytical solution of the DRE developed by Liao \cite{liao2021research}, the analytical strategy presented in this paper is more generalizable, since the matrix $\bm{S}$ is not required to be strictly positive definite, and the developed strategy is applicable to arbitrary Keplerian reference orbits.

\noindent \textbf{Remark 2}: Compared to the game strategy based on the numerical solution of the DRE \cite{li2024nash} (i.e., the numerical strategy), the analytical strategy developed in this paper offers higher computational efficiency.

\noindent \textbf{Remark 3}: While the analytical strategy is developed for the pursuit-evasion game scenario in this paper, a similar analytical approach can be adopted to solve the DRE for other scenarios, such as terminal guidance for NEAs and space debris defense. For such extensions, some modifications are necessary, which are detailed in Appendix \ref{secA4}.

Subsequently, a procedure for solving the LQ pursuit-evasion game using the analytical strategy is proposed.

\subsection{Procedure Using the Analytical Strategy}\label{subsec3.2}
\begin{figure}[H]
\centering
\includegraphics[width=0.8\textwidth]{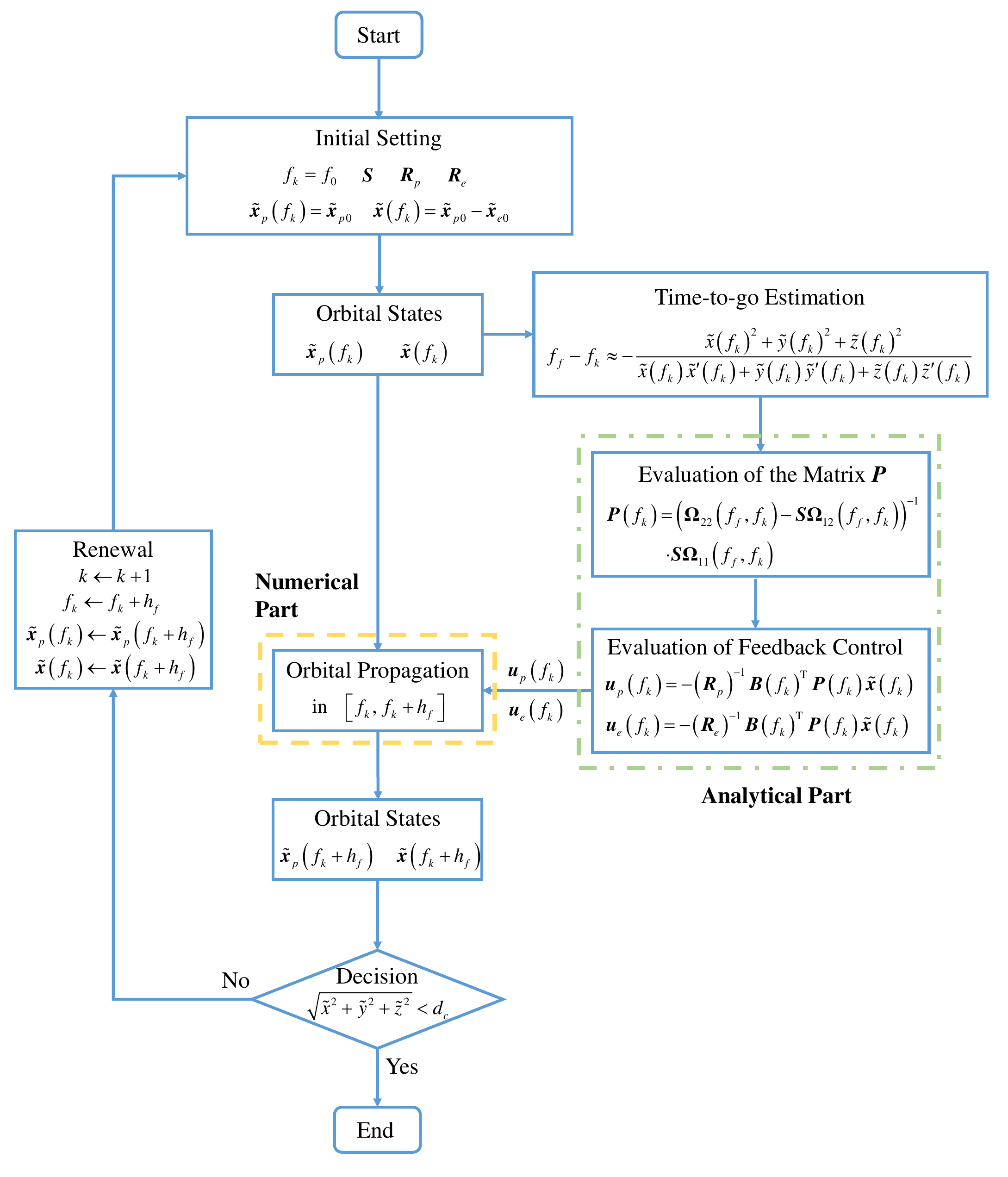}
\caption{Flow Chart of the Procedure Using the Analytical Strategy.}\label{fig4}
\end{figure}

The flow chart of the procedure using the analytical strategy is presented in Fig. \ref{fig4}. In this procedure, when the relative states between the pursuer and evader satisfy $\sqrt {{{\tilde x}^2} + {{\tilde y}^2} + {{\tilde z}^2}}  < {d_c}$, the evader is captured by the pursuer. In the capture decision, $d_c$ represents the capture radius for the pursuit-evasion game. Since the time-to-go estimation and the capture decision are required when solving the pursuit-evasion game, the dynamical equations Eq. \eqref{eq7}-\eqref{eq8} is solved numerically although $\bm{u}_p$ and $\bm{u}_e$ can be calculated analytically once the $\tilde{\bm x}$ is obtained.

\subsection{Limitations}\label{subsec3.3}
When the analytical strategy and the procedure for solving the LQ pursuit-evasion game are proposed, several potential limitations should be acknowledged:

\begin{enumerate}[label=(\arabic*)]
\item Since the analytical strategy is presented based on linear relative dynamics (i.e., controlled TH equations), it may have limitations in scenarios involving strong nonlinear perturbations, such as multi-body perturbation.
\item The analytical strategy is developed under the assumption of complete information. In practical applications, however, sensor noise or measurement delays could affect access to accurate state information, thereby affecting the accuracy of the control input. Small errors can accumulate over time, which may be particularly problematic in long-duration pursuit-evasion scenarios.
\end{enumerate}

\section{Results and Discussion}\label{sec4}
In this section, simulations of the pursuit-evasion game in elliptic, parabolic, and hyperbolic reference orbits are performed using the developed analytical strategy. The results obtained from the analytical strategy are compared with those obtained from the numerical strategy to illustrate the efficiency of the developed analytical strategy. Furthermore, the effectiveness of the developed analytical strategy in the pursuit-evasion game considering orbital disturbances is analyzed. Finally, a comparison between the developed and conventional strategies is presented to highlight the advantages of the developed strategy.
\subsection{Parameter Setting for Simulation}\label{subsec4.1}
The parameters used to solve the pursuit-evasion game in the simulations are specified in Table \ref{tab1}. The CPU used for the simulations is the Intel® i9-13900KF processor with a base frequency of 3.00 GHz and 64 GB of RAM. The software employed in the simulations is Matlab®. To demonstrate the performance of the proposed strategy, two pursuit-evasion game scenarios are discussed in this paper, denoted as Scenario I and Scenario II. The parameters of these scenarios are summarized in Tables \ref{tab2} and \ref{tab3}. Simulations for these two scenarios are then performed in three cases, corresponding to elliptic, parabolic, and hyperbolic reference orbits, respectively.

\begin{table}[!htb]
\caption{Parameter setting for solving the pursuit-evasion game.}\label{tab1}%
\centering
\renewcommand{\arraystretch}{1.5}
\begin{tabular}{@{}llll@{}}
\hline
Symbol & Value  & Units & Meaning\\
\hline
$\mu$    & $3.98603 \times {10^{14}}$   & $\text{m}{}^3/\text{s}{}^2$  & Earth Gravitational Constant  \\
$d_c$ & $1$ & $\text{m}$ & Capture Radius for the Pursuit-Evasion Game  \\
\hline
\end{tabular}
\end{table}

\begin{table}[!htb]
\caption{Parameter setting for Scenario I.}\label{tab2}%
\centering
\renewcommand{\arraystretch}{1.5}
\begin{tabular}{@{}llll@{}}
\hline
Symbol & Value  & Units & Meaning\\
\hline
$p$   & $4.2241 \times {10^{7}}$   & $\text{m}$  & Semilatus Rectum of the reference orbit  \\
$s_r$    & $0.1$   & --  & Weight of the Weighting Matrix $\bm{S}$  \\
$r_p$    & $1 \times {10^{6}}$   & --  & Weight of the Weighting Matrix $\bm{R}_p$  \\
$r_e$    & $1.1 \times {10^{6}}$  & --  & Weight of the Weighting Matrix $\bm{R}_e$  \\
$[\tilde{x}_{p0},\text{ }\tilde{y}_{p0},\text{ }\tilde{z}_{p0}]$    & $[1500,\text{ }500,\text{ }0]$   & $\text{m}$  & Initial Position of Pursuer in the LVLH Frame  \\
$[{\tilde x_{p0}}',\text{ }{\tilde y_{p0}}',\text{ }{\tilde z_{p0}}'] $    & $[-10000,\text{ }0,\text{ }1000]$   & \text{m/rad}  &Initial Velocity of Pursuer in the LVLH Frame  \\
$[\tilde{x}_{e0},\text{ }\tilde{y}_{e0},\text{ }\tilde{z}_{e0}]$    & $[0,\text{ }0,\text{ }0]$   & $\text{m}$  & Initial Position of Evader in the LVLH Frame  \\
$[{\tilde x_{e0}}',\text{ }{\tilde y_{e0}}',\text{ }{\tilde z_{e0}}'] $    & $[0,\text{ }0,\text{ }0]$   & m/rad  &Initial Velocity of Evader in the LVLH Frame  \\
$f_0$ & 0   & $\text{rad}$  & Initial True Anomaly of the Pursuit-Evasion Game  \\
\hline
\end{tabular}
\end{table}

\begin{table}[!htb]
\caption{Parameter setting for Scenario II.}\label{tab3}%
\centering
\renewcommand{\arraystretch}{1.5}
\begin{tabular}{@{}llll@{}}
\hline
Symbol & Value  & Units & Meaning\\
\hline
$p$   & $1.8339 \times {10^{7}}$   & $\text{m}$  & Semilatus Rectum of the reference orbit  \\
$s_r$    & $1$   & --  & Weight of the Weighting Matrix $\bm{S}$  \\
$r_p$    & $1 \times {10^{6}}$   & --  & Weight of the Weighting Matrix $\bm{R}_p$  \\
$r_e$    & $1.5 \times {10^{6}}$  & --  & Weight of the Weighting Matrix $\bm{R}_e$  \\
$[\tilde{x}_{p0},\text{ }\tilde{y}_{p0},\text{ }\tilde{z}_{p0}]$   & $[1500,\text{ }0,\text{ }-2500]$   & $\text{m}$  & Initial Position of Pursuer in the LVLH Frame  \\
$[{\tilde x_{p0}}',\text{ }{\tilde y_{p0}}',\text{ }{\tilde z_{p0}}'] $    & $[-5000,\text{ }0,\text{ }10000]$   & \text{m/rad}  &Initial Velocity of Pursuer in the LVLH Frame  \\
$[\tilde{x}_{e0},\text{ }\tilde{y}_{e0},\text{ }\tilde{z}_{e0}]$    & $[0,\text{ }0,\text{ }0]$   & $\text{m}$  & Initial Position of Evader in the LVLH Frame  \\
$[{\tilde x_{e0}}',\text{ }{\tilde y_{e0}}',\text{ }{\tilde z_{e0}}'] $    & $[0,\text{ }0,\text{ }0]$   & m/rad  &Initial Velocity of Evader in the LVLH Frame \\
$f_0$ & 0   & $\text{rad}$  & Initial True Anomaly of the Pursuit-Evasion Game  \\
\hline
\end{tabular}
\end{table}

\subsection{Effectiveness of the Procedure Using the Analytical Strategy}\label{subsec4.2}

The procedure using the analytical strategy is shown in Fig. \ref{fig4}. The dynamical equations Eqs. \eqref{eq7}-\eqref{eq8} are solved by the fourth-order Runge-Kutta method \cite{battin1999an} with a step size of $h_f=0.00001$. To highlight the effectiveness and efficiency of the analytical strategy, a comparison  is presented between the results obtained from the analytical and numerical strategies. Firstly, the procedure using the numerical strategy is proposed. In this procedure, the control inputs of the pursuer and evader $\bm{u}_p$ and $\bm{u}_e$ at $f$ are calculated using the numerical solution of the DRE Eq. \eqref{eq22}. The DRE is solved numerically with the same time-to-go estimation mentioned in Section \ref{subsubsec3.1.4}. For the procedure using the numerical strategy, the fourth-order Runge-Kutta method is adopted to solve both the dynamical equations Eqs. \eqref{eq7}-\eqref{eq8} and the DRE Eq. \eqref{eq22} (note that the step size to solve the dynamical equations is the same as that used in the procedure using the analytical strategy). Subsequently, the simulations are performed for the aforementioned two scenarios in elliptic, parabolic, and hyperbolic reference orbits, respectively.

\subsubsection{Case I}\label{subsubsec4.2.1}

\begin{figure}[H]
\centering
\includegraphics[width=0.64\textwidth]{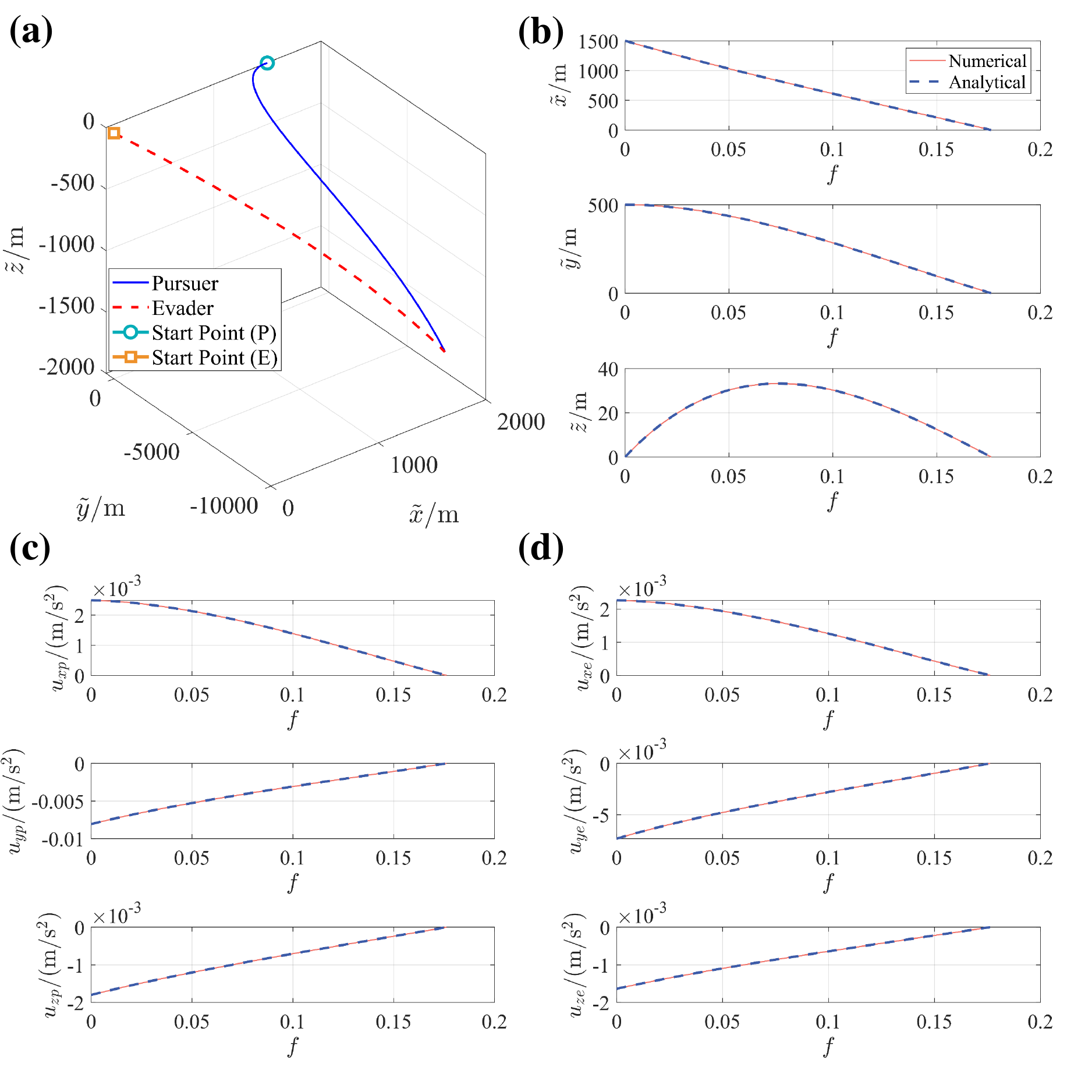}
\caption{Pursuit-Evasion game in an elliptic reference orbit (Scenario I). (a) Trajectories of the pursuer and evader. (b) Relative position. (c) Control input of the pursuer. (d) Control input of the evader.}\label{fig_elliptic_sc_I}
\end{figure}

\begin{figure}[H]
\centering
\includegraphics[width=0.64\textwidth]{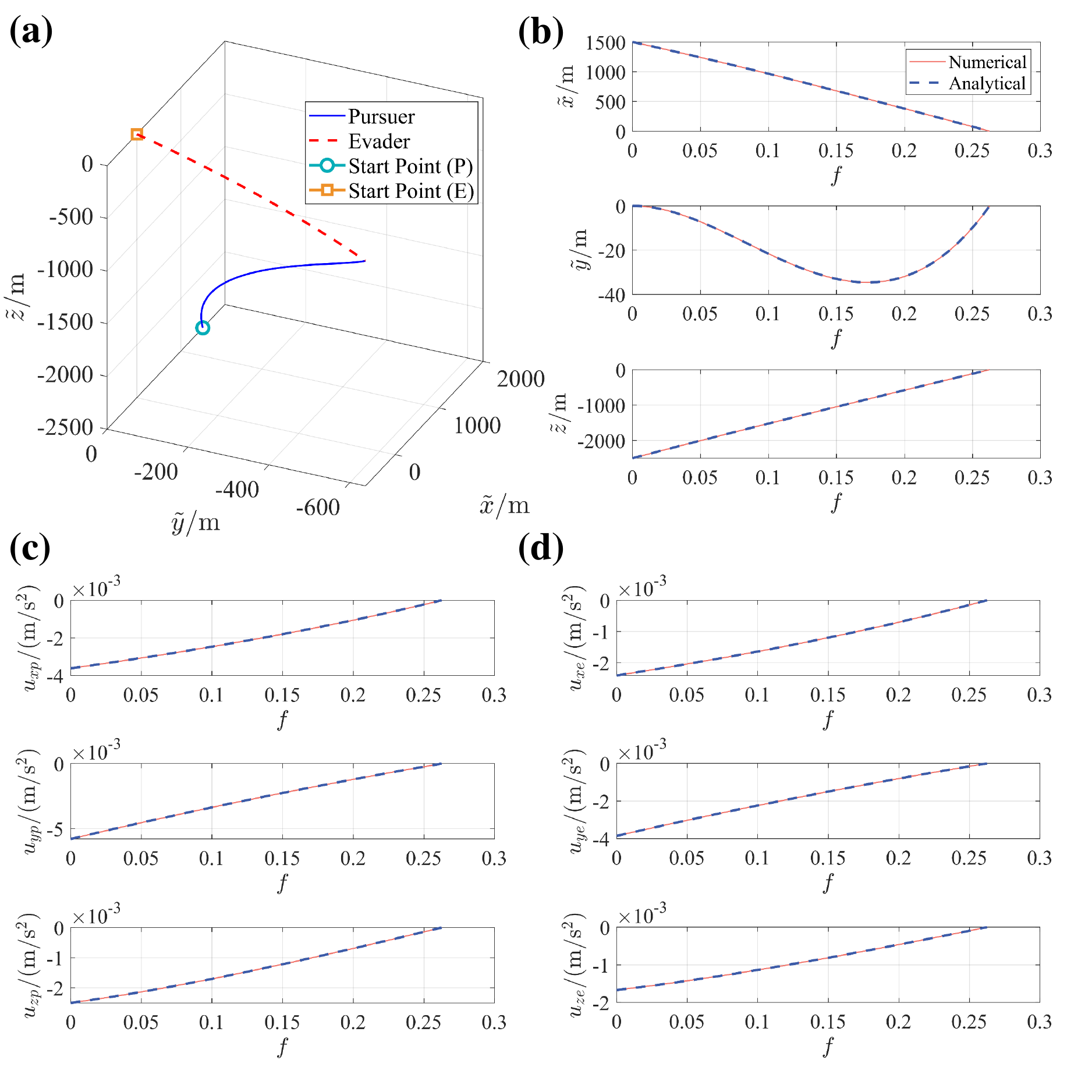}
\caption{Pursuit-Evasion game in an elliptic reference orbit (Scenario II). (a) Trajectories of the pursuer and evader. (b) Relative position. (c) Control input of the pursuer. (d) Control input of the evader.}\label{fig_elliptic_sc_II}
\end{figure}

For the pursuit-evasion game in an elliptic reference orbit, we set the eccentricity of the reference orbit to $0.2$. The developed analytical strategy is also applicable to the circular reference orbit ($e=0$). The trajectories of the pursuer and evader obtained from the analytical strategy are shown in Figs. \ref{fig_elliptic_sc_I} and \ref{fig_elliptic_sc_II} (a). It is found that both the pursuer and evader move away from the reference orbit, and the evader tends to move in the $\tilde{y}$ direction relative to the pursuer in these two scenarios. Finally, the evader is captured by the pursuer. Figures \ref{fig_elliptic_sc_I} and \ref{fig_elliptic_sc_II} (b) present the relative position components between the pursuer and evader during the game obtained from the analytical and numerical strategies. It is observed that the results obtained from the analytical strategy are in good agreement with those obtained from the numerical one. This consistency is also evident in the control inputs of the pursuer and evader shown in Figs. \ref{fig_elliptic_sc_I} and \ref{fig_elliptic_sc_II} (c) and (d). Therefore, the effectiveness of the analytical strategy is validated. Furthermore, for Scenario I, the CPU time for the procedure using the numerical strategy is 2117.45 s, while the CPU time for the procedure using the analytical strategy is only 3.63 s. For Scenario II, the CPU time for the procedure using the numerical strategy is 3166.25 s, while the procedure using the analytical strategy requires just 5.68 s. Therefore, it can be concluded that the analytical strategy saves over 99.82 $\%$ of the CPU time compared to the numerical one for the two scenarios in elliptic reference orbits, while maintaining control accuracy. It can be concluded that the numerical integration of the DRE Eq. \eqref{eq22} is the main factor contributing to the high CPU time in the procedure using the numerical strategy.

\subsubsection{Case II}\label{subsec4.2.2}
\begin{figure}[H]
\centering
\includegraphics[width=0.64\textwidth]{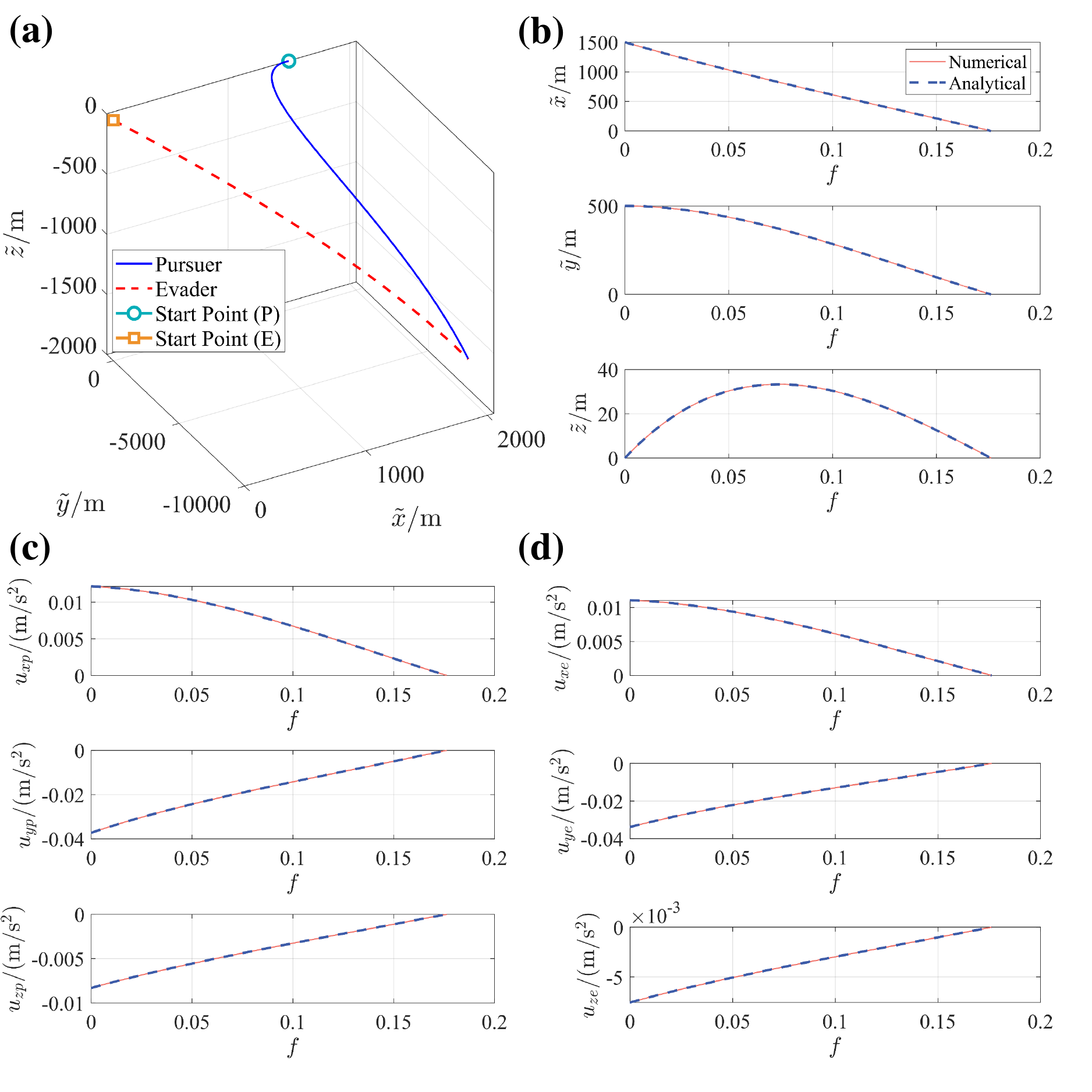}
\caption{Pursuit-Evasion game in a parabolic reference orbit (Scenario I). (a) Trajectories of the pursuer and evader. (b) Relative position. (c) Control input of the pursuer. (d) Control input of the evader.}\label{fig_parabolic_sc_I}
\end{figure}

\begin{figure}[H]
\centering
\includegraphics[width=0.64\textwidth]{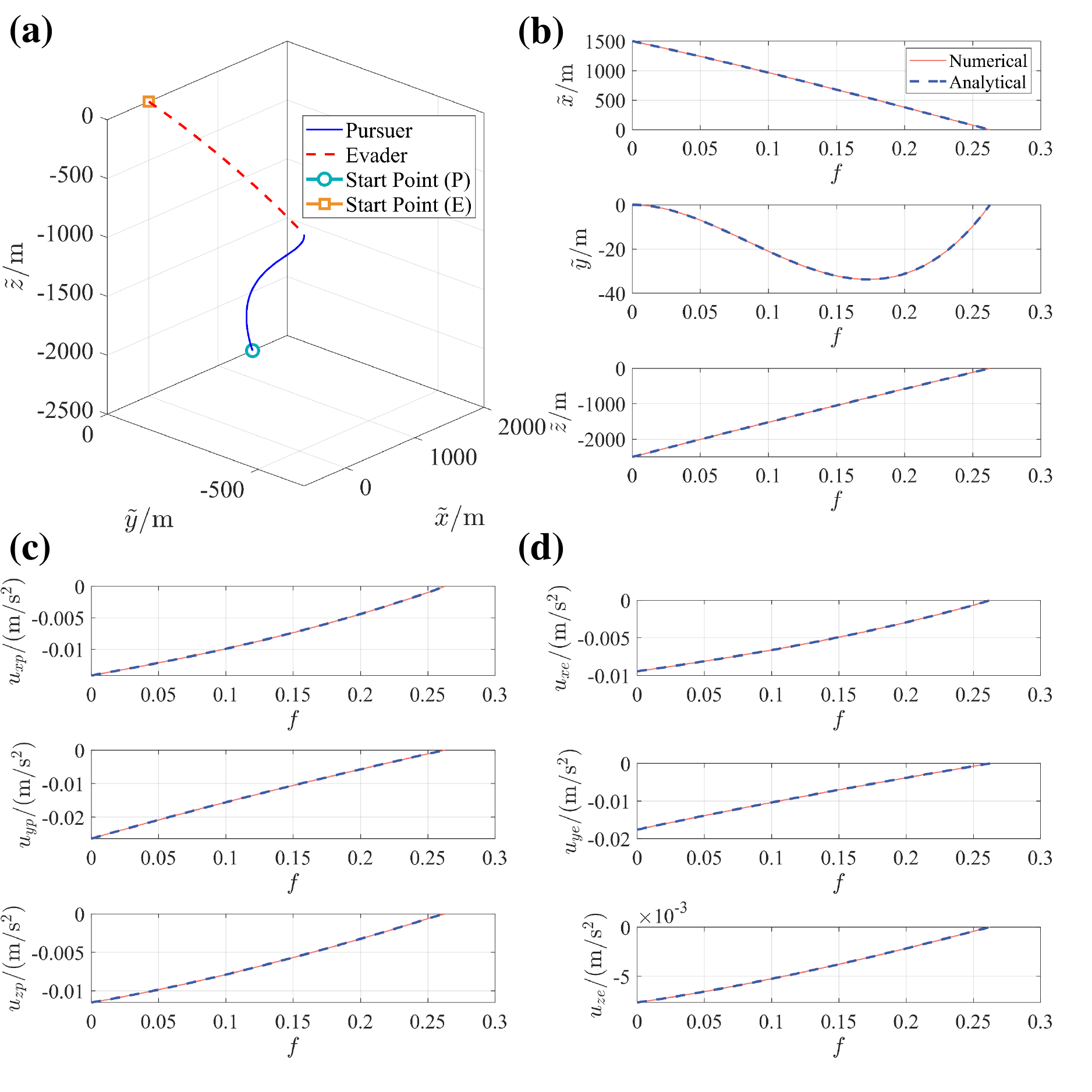}
\caption{Pursuit-Evasion game in a parabolic reference orbit (Scenario II). (a) Trajectories of the pursuer and evader. (b) Relative position. (c) Control input of the pursuer. (d) Control input of the evader.}\label{fig_parabolic_sc_II}
\end{figure}

When the reference orbit is a parabolic orbit, the evader is also captured by the pursuer in Scenarios I and II as in the elliptic reference orbit, as shown in Figs. \ref{fig_parabolic_sc_I} and \ref{fig_parabolic_sc_II} (a) and (b). Differences are observed in the control inputs of the pursuer and evader due to the high eccentricity and increased speed near the reference orbit, as presented in Figs. \ref{fig_parabolic_sc_I} and \ref{fig_parabolic_sc_II}(c) and (d). A good agreement between the results obtained from the analytical and numerical strategies is evident, demonstrating the accuracy of the proposed procedure. For Scenario I, the CPU time required by the procedure using the numerical strategy is 2121.94 s, while the procedure using the analytical strategy requires only 3.09 s. Similarly, for Scenario II, the CPU time for the procedure using the numerical strategy is 3164.39 s, compared to just 4.78 s for the procedure using the analytical strategy. In the parabolic reference orbit, the analytical strategy saves the computation time by more than 99.84$\%$ compared to the numerical one.

\subsubsection{Case III}\label{subsec4.2.3}
\begin{figure}[H]
\centering
\includegraphics[width=0.64\textwidth]{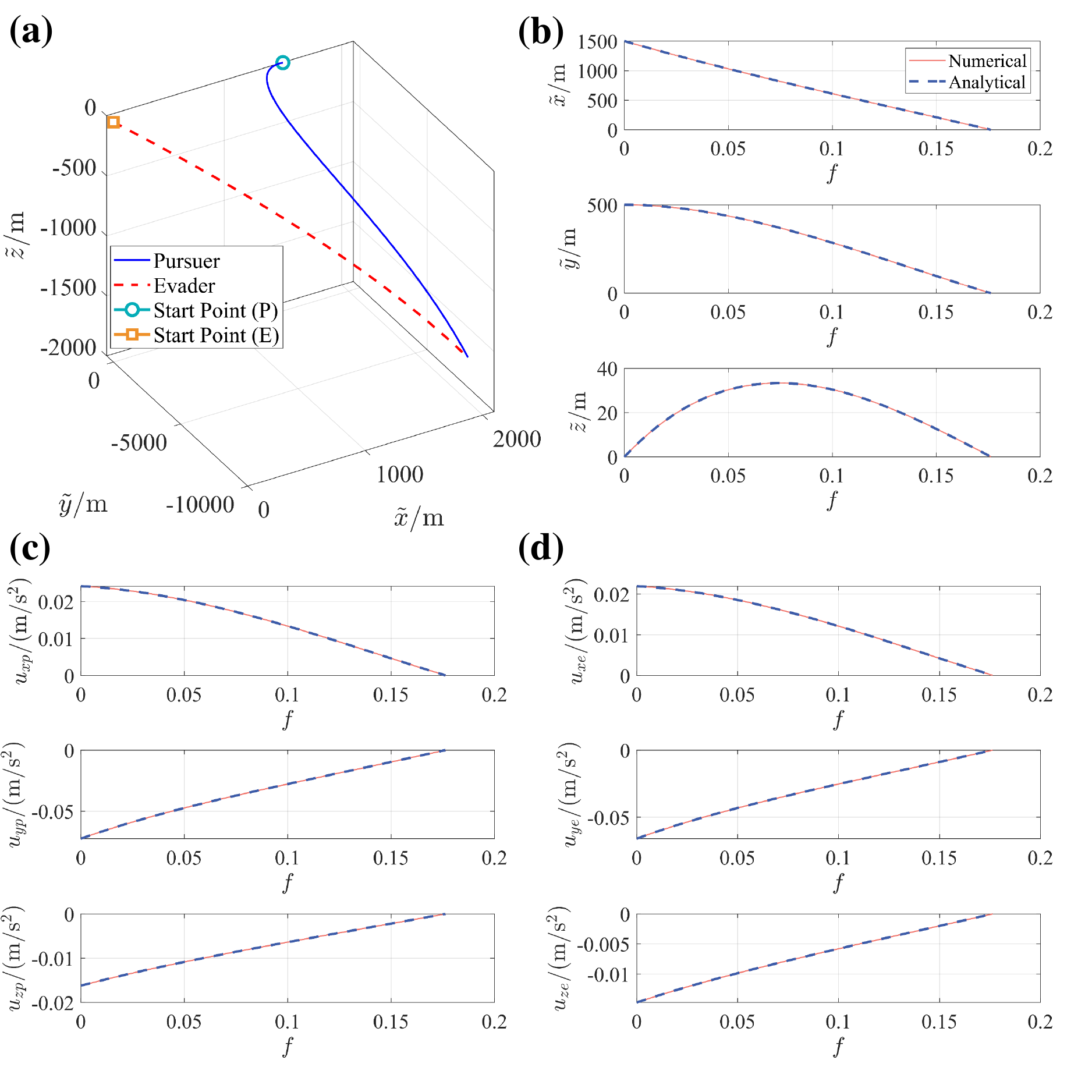}
\caption{Pursuit-Evasion game in a hyperbolic reference orbit (Scenario I). (a) Trajectories of the pursuer and evader. (b) Relative position. (c) Control input of the pursuer. (d) Control input of the evader.}\label{fig_hyperbolic_sc_I}
\end{figure}

\begin{figure}[H]
\centering
\includegraphics[width=0.64\textwidth]{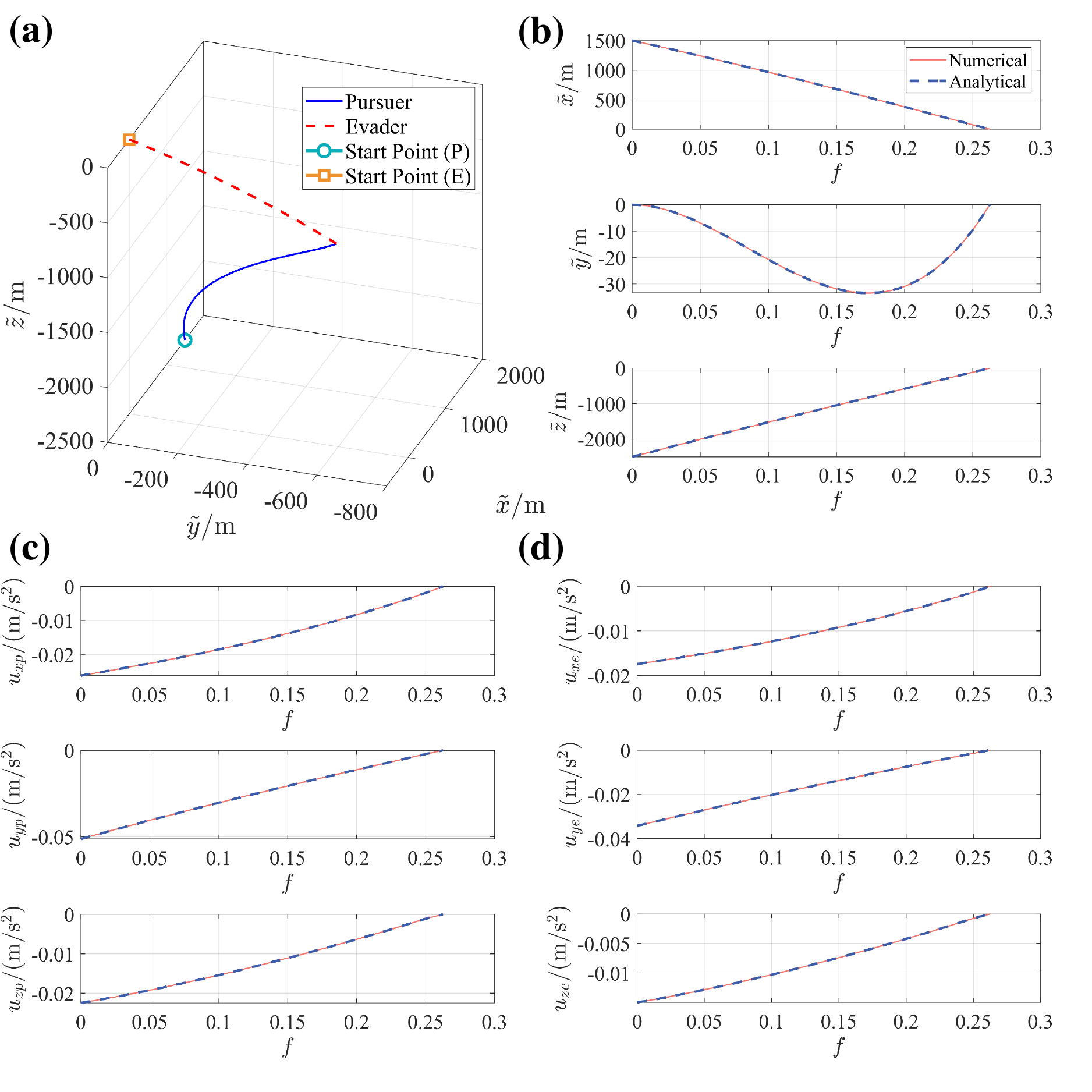}
\caption{Pursuit-Evasion game in a hyperbolic reference orbit (Scenario II). (a) Trajectories of the pursuer and evader. (b) Relative position. (c) Control input of the pursuer. (d) Control input of the evader.}\label{fig_hyperbolic_sc_II}
\end{figure}

For the case of a hyperbolic reference orbit, the eccentricity is set to 1.5. Figures \ref{fig_hyperbolic_sc_I} and \ref{fig_hyperbolic_sc_II} (a) show the trajectories of the pursuer and evader, and Figs. \ref{fig_hyperbolic_sc_I} and \ref{fig_hyperbolic_sc_II} (b)-(d) present a comparison of their relative position components and control inputs. A good agreement between the results obtained from the analytical and numerical strategies is observed. In the hyperbolic reference orbit, the CPU time required by the procedure using the numerical strategy is 2108.85 s, while the procedure using the analytical strategy takes only 3.69 s for Scenario I. For Scenario II, it is 3151.33 s for the procedure using the numerical strategy and 5.81 s for the procedure using the analytical strategy. For the scenarios in the hyperbolic reference orbits, the analytical strategy saves the CPU time by more than 99.81$\%$ compared to the numerical one.

Above all, to quantify the differences between the results obtained from the analytical and the numerical strategies, three key metrics are recorded for the two scenarios in the elliptic, parabolic, and hyperbolic reference orbits (i.e., six test cases). These metrics include the position error at the terminal true anomaly $f_f$ (i.e., $\tilde r\left( {{f_f}} \right)=\sqrt {{{\tilde x}^2}\left( {{f_f}} \right) + {{\tilde y}^2}\left( {{f_f}} \right) + {{\tilde z}^2}\left( {{f_f}} \right)} $), the game time (i.e., $\Delta{f}={f_f}-{f_0}$), and the cost function (i.e., $J$ defined in Eq. \eqref{eq9}). The detailed results for these metrics are summarized in Table \ref{tab4}.

\begin{table}[!htb]
\renewcommand{\arraystretch}{1.5}
\caption{Quantitative differences between the results obtained from the analytical and numerical strategies.}\label{tab4}%
	\centering

	\begin{tabular}{@{}llllll@{}}
 \hline
		
		Strategy     &         Scenario, Type of Reference Orbit                      & $\tilde r\left( {{f_f}} \right)/\text{m} $    & $\Delta{f}$ & $J$ & CPU Time/s \\ \hline
		\multirow{6}{*}{Numerical} & \multicolumn{1}{l}{Scenario I, Elliptic} & 0.9914 & 0.17621 & 0.2298 & 2117.45   \\
		                     & \multicolumn{1}{l}{Scenario I, Parabolic} & 0.9541 & 0.17625 & 3.9538 & 2121.94   \\
		                     & \multicolumn{1}{l}{Scenario I, Hyperbolic} & 0.9790 & 0.17627 & 14.9915 & 2108.85   \\
		                     & \multicolumn{1}{l}{Scenario II, Elliptic} & 0.9689 & 0.26255 & 1.2399 & 3166.25   \\ 
                              & \multicolumn{1}{l}{Scenario II, Parabolic} & 0.9380 & 0.26259 & 15.6849 & 3164.39   \\
                             & \multicolumn{1}{l}{Scenario II, Hyperbolic} & 0.9288 & 0.26259 & 57.2257 & 3151.33   \\ \hline
		\multirow{6}{*}{Analytical}& \multicolumn{1}{l}{Scenario I, Elliptic} & 0.9755 & 0.17615 & 0.2282 & 3.63   \\
		                     & \multicolumn{1}{l}{Scenario I, Parabolic} & 0.9820 & 0.17620 & 3.9560 & 3.09  \\
		                     & \multicolumn{1}{l}{Scenario I, Hyperbolic} & 0.9772 & 0.17623 & 14.9894 & 3.69   \\
		                     & \multicolumn{1}{l}{Scenario II, Elliptic} & 0.9731 & 0.26249 & 1.2439 & 5.68   \\ 
                              & \multicolumn{1}{l}{Scenario II, Parabolic} & 0.9990 & 0.26254 & 15.7472 & 4.78   \\
                             & \multicolumn{1}{l}{Scenario II, Hyperbolic} & 0.9520 & 0.26255 & 57.2428 & 5.81   \\ \hline
		
	\end{tabular}
	
\end{table}

Across these six test cases, the analytical strategy saves the CPU time by more than 99.80$\%$ compared to the numerical one. Meanwhile, as shown in Table \ref{tab4}, there are no significant differences in the position error, the game time, and the cost function between the results obtained from the analytical and numerical strategies. This consistency validates the effectiveness of the analytical strategy while highlighting its advantage in computational efficiency.

\subsection{Cases with Orbital Disturbances}\label{subsec4.4}
\begin{figure}[H]
\centering
\includegraphics[width=0.64\textwidth]{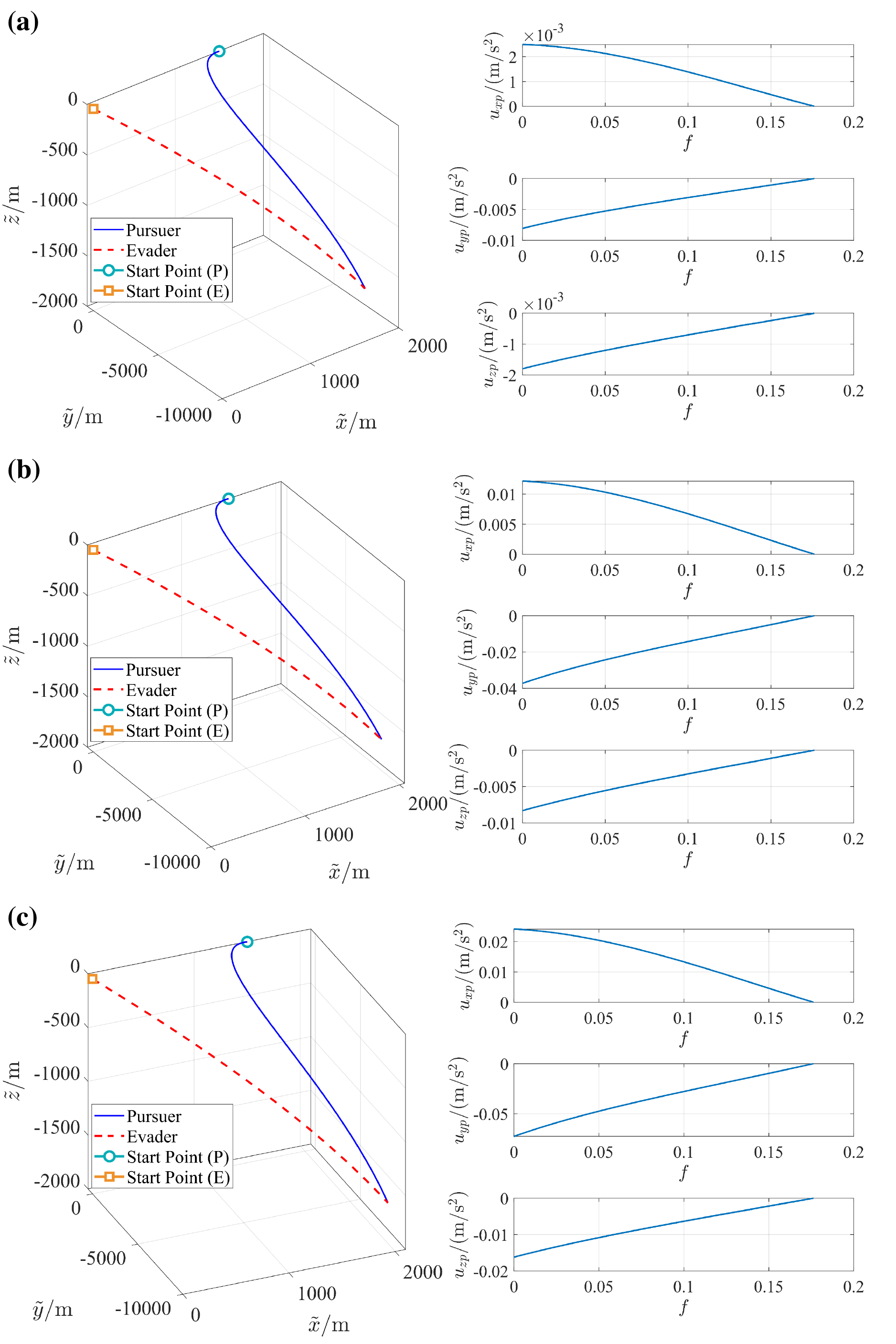}
\caption{Trajectories of the pursuer and evader with the orbital disturbances and the variation in the control input of the pursuer (Scenario I). (a) Elliptic reference orbit. (b) Parabolic reference orbit. (c) Hyperbolic reference orbit.}\label{sc_I_disturbances}
\end{figure}

\begin{figure}[H]
\centering
\includegraphics[width=0.64\textwidth]{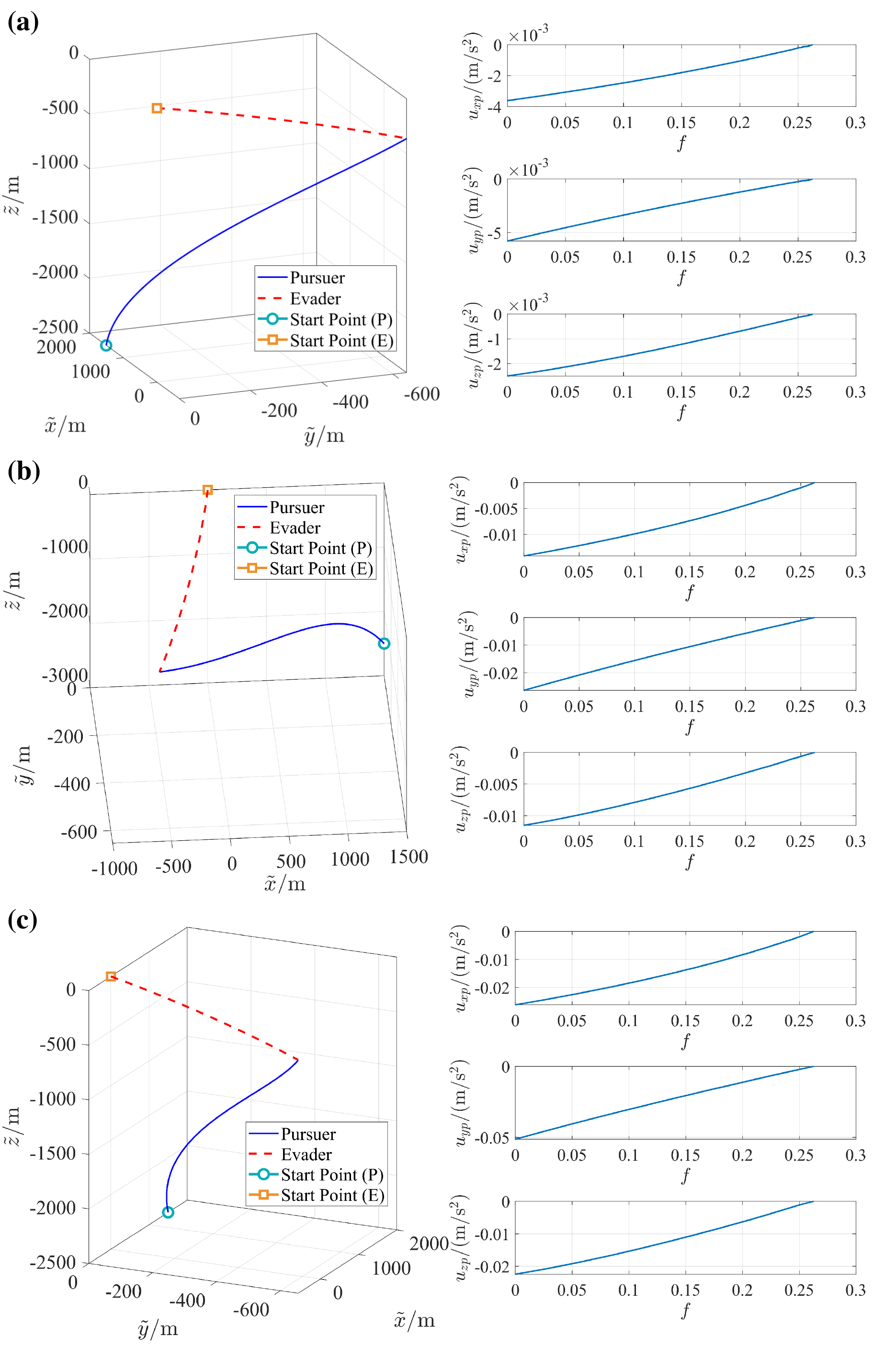}
\caption{Trajectories of the pursuer and evader with the orbital disturbances and the variation in the control input of the pursuer (Scenario II). (a) Elliptic reference orbit. (b) Parabolic reference orbit. (c) Hyperbolic reference orbit.}\label{sc_II_disturbances}
\end{figure}

Based on the aforementioned discussion, it is concluded that the analytical strategy performs well in all six test cases. Therefore, in Section \ref{subsec4.4}, the discussion continues with these test cases. To further evaluate the performance of the developed analytical strategy, orbital disturbances are introduced into the TH equations of the pursuer and evader:

\begin{equation}
\left\{ \begin{gathered}
  {{\tilde{\bm X}'}_p} = \bm{A}{{\tilde{\bm X}}_p} + \bm{B}{\bm{u}_p} +{\bm{d}_p}  \hfill \\
  {{\tilde{\bm X}'}_e} = \bm{A}{{\tilde{\bm X}}_e} + \bm{B}{\bm{u}_e} +{\bm{d}_e} \hfill \\ 
\end{gathered}  \right.\label{eq10000}
\end{equation}
where $\bm{d}_p=[0,\text{ }0,\text{ }0,\text{ }d_{px},\text{ }d_{py},\text{ }d_{pz}]^{\text{T}}$ and $\bm{d}_e=[0,\text{ }0,\text{ }0,\text{ }d_{ex},\text{ }d_{ey},\text{ }d_{ez}]^{\text{T}}$ denote the orbital disturbances of the pursuer and evader. Then, the dynamical equations between the purser and evader considering orbital disturbances can be expressed as:

\begin{equation}
\tilde{\bm x}' = \bm{A}\tilde{\bm x} + \bm{B}{\bm{u}_p} - \bm{B}{\bm{u}_e} + {\bm{d}_p}- {\bm{d}_e}\label{eq10001}
\end{equation}
In the pursuit-evasion game considering orbital disturbances, the control inputs $\bm{u}_p$ and $\bm{u}_e$ are still calculated analytically by the Eq. \eqref{eq12} when $\tilde{\bm{x}}$ is provided. The procedure presented in Fig. \ref{fig4} (analytically obtaining the control inputs and numerically solving the Eqs. \eqref{eq10000}-\eqref{eq10001} under orbital disturbances) is performed. The disturbances are randomly generated, where $d_{ab}\text{ }(a=p,\text{ }e;\text{ }b=x,\text{ }y,\text{ }z)$ satisfies $d_{ab} \in \left[-1000,\text{ }1000\right]\text{ }\left(\text{m/rad}^2\right)$, achieved by Matlab® ‘rand’ command. For the six test cases, the evader is successfully captured by the pursuer, as shown in Figs. \ref{sc_I_disturbances}-\ref{sc_II_disturbances}. The variations in $\bm{u}_p$ in these six test cases are also presented (according to Eq. \eqref{eq12}, the control inputs of the pursuer and evader are proportional, i.e., $\bm{u}_p=\alpha\bm{u}_e$; therefore, only the variations in $\bm{u}_p$ are presented). It can be concluded that the analytical strategy remains applicable even in scenarios under orbital disturbances.

\subsection{Comparison with Conventional Strategy}\label{subsec4.5}
\begin{figure}[H]
\centering
\includegraphics[width=0.64\textwidth]{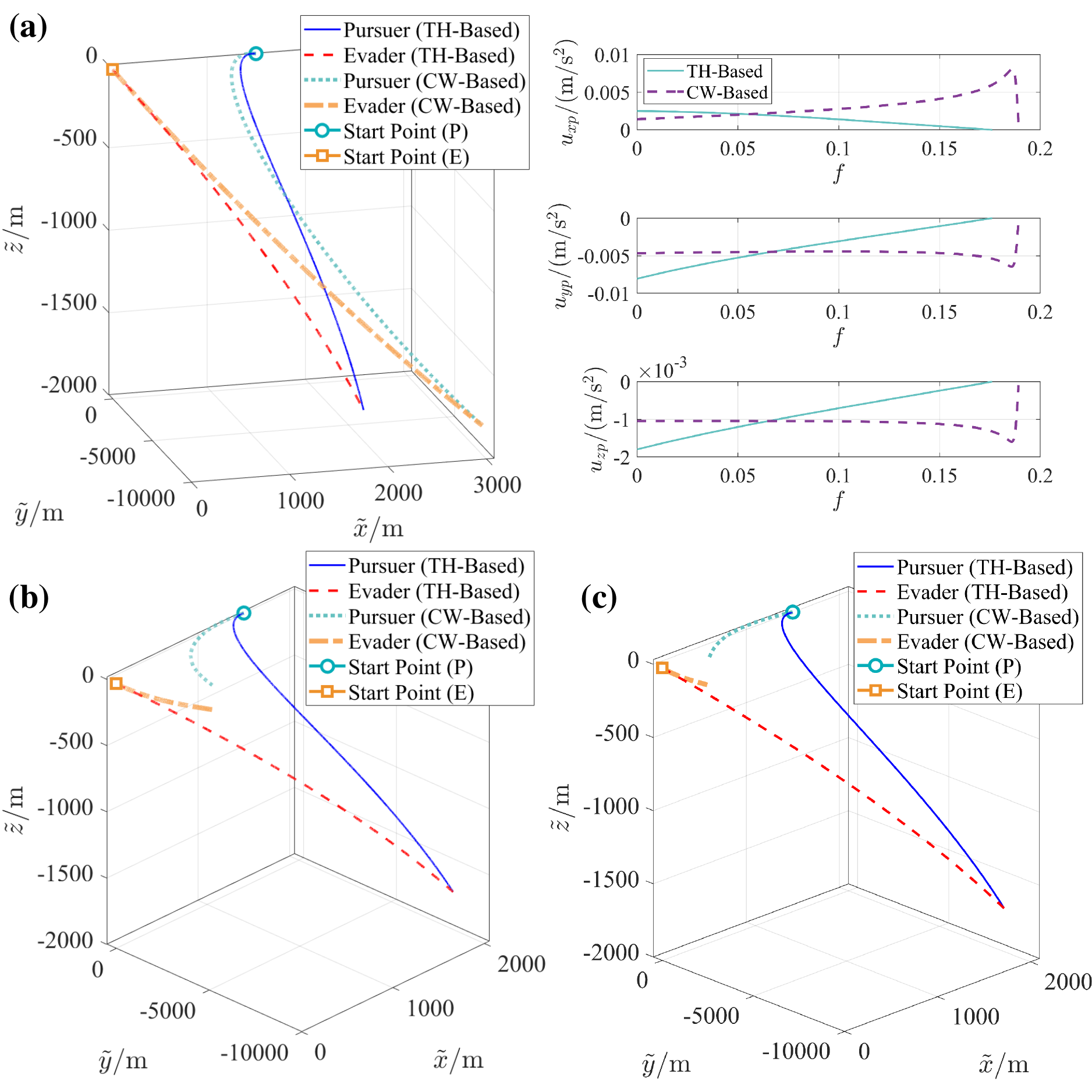}
\caption{Comparison between the TH-based strategy and CW-based strategy (Scenario I). (a) Trajectories in the elliptic reference orbit and the variation in $\bm{u}_p$. (b) Trajectories in the parabolic reference orbit. (c) Trajectories in the hyperbolic reference orbit.}\label{Fig451}
\end{figure}

\begin{figure}[H]
\centering
\includegraphics[width=0.64\textwidth]{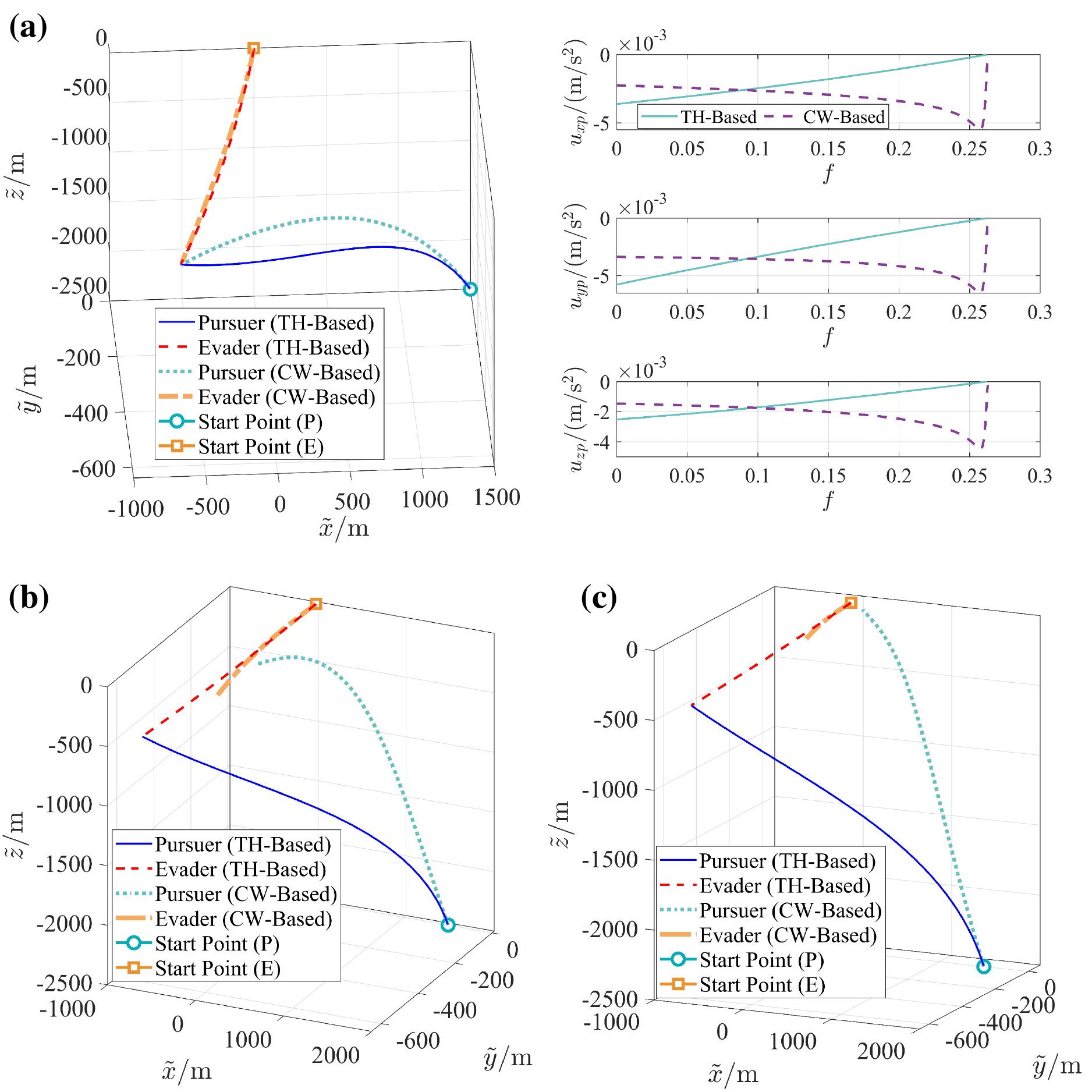}
\caption{Comparison between the TH-based and CW-based strategy (Scenario II). (a) Trajectories in the elliptic reference orbit and the variation in $\bm{u}_p$. (b) Trajectories in the parabolic reference orbit. (c) Trajectories in the hyperbolic reference orbit.}\label{Fig452}
\end{figure}

When the effectiveness of the developed startegy in the pursuit-evasion game considering orbital disturbances is discussed, a comparison with the conventional strategies can further highlight the advantage of the strategy developed in this paper. The aforementioned six test cases (i.e., Scenarios I and II in the elliptic ($e=0.2$), parabolic ($e=1$), and hyperbolic ($e=1.5$) reference orbits) are used for this comparison. The conventional strategy considered here is the CW-based strategy \cite{liao2021research} (i.e., in this strategy, the control input is generated from the analytical solution of the DRE Eq. \eqref{eq22} in the circular reference orbit). In contrast, the strategy developed in this paper is denoted as the TH-based strategy. The results of this comparison are presented in Figs. \ref{Fig451}-\ref{Fig452}. For the pursuit-evasion game in elliptic reference orbits, the evader can be successfully captured by the pursuer using the CW-based strategy. However, the CW-based strategy fails in both parabolic and hyperbolic reference orbits for Scenarios I and II, whereas the TH-based strategy remains effective across all test cases.

To further highlight the advantage of the TH-based strategy, the position error, the game time and the cost function are also recorded and summarized in Table \ref{tab5}. For the failed test cases, the position error represents the minimum relative distance between the pursuer and evader during the game.

\begin{table}[!htb]
\renewcommand{\arraystretch}{1.5}
\caption{Comparison between the TH-based and CW-based strategies.}\label{tab5}%
	\centering

	\begin{tabular}{@{}lllll@{}}
 \hline
		
		Strategy     &         Scenario, Type of Reference Orbit                      & $\tilde r\left( {{f_f}} \right)/\text{m} $    & $\Delta{f}$ & $J$  \\ \hline
		\multirow{6}{*}{TH-based}& \multicolumn{1}{l}{Scenario I, Elliptic} & 0.9755 & 0.17615 & 0.2282   \\
		                     & \multicolumn{1}{l}{Scenario I, Parabolic} & 0.9820 & 0.17620 & 3.9560   \\
		                     & \multicolumn{1}{l}{Scenario I, Hyperbolic} & 0.9772 & 0.17623 & 14.9894   \\
		                     & \multicolumn{1}{l}{Scenario II, Elliptic} & 0.9731 & 0.26249 & 1.2439  \\ 
                              & \multicolumn{1}{l}{Scenario II, Parabolic} & 0.9990 & 0.26254 & 15.7472   \\
                             & \multicolumn{1}{l}{Scenario II, Hyperbolic} & 0.9520 & 0.26255 & 57.2428    \\ \hline
        \multirow{6}{*}{CW-based}& \multicolumn{1}{l}{Scenario I, Elliptic} & 0.9574 & 0.18914 & 0.3442    \\
		                     & \multicolumn{1}{l}{Scenario I, Parabolic} & 589.9189 & -- & --   \\
		                     & \multicolumn{1}{l}{Scenario I, Hyperbolic} & 656.1038 & -- & --   \\
		                     & \multicolumn{1}{l}{Scenario II, Elliptic} & 0.9253 & 0.26270 & 1.7730  \\ 
                              & \multicolumn{1}{l}{Scenario II, Parabolic} & 212.6866 & -- & --   \\
                             & \multicolumn{1}{l}{Scenario II, Hyperbolic} & 329.2829 & -- & --    \\ \hline
		
	\end{tabular}
	
\end{table}

From the table, it is observed that, for the successful test cases, the game times obtained using the CW-based strategy are longer than those obtained using the TH-based strategy, and the values of the cost function are also higher when using the CW-based strategy. For the failure test cases, it is observed that as the value of the eccentricity increases, the minimum relative distance using the CW-based strategy also increases for both Scenarios I and II. Above all, the CW-based strategy succeeds in only two out of the six test cases while the TH-based strategy achieves success in all six test cases. This demonstrates that the TH-based strategy developed in this paper is more effective for solving the pursuit-evasion game in high-eccentricity reference orbits.

\section{Conclusion}\label{sec5}
 
This paper extends the linear quadratic (LQ) pursuit-evasion game to arbitrary Keplerian reference orbits, including elliptic, parabolic, and hyperbolic orbits. For each case, the analytical strategy for the pursuit-evasion game is developed based on the analytical solution of the differential Riccati equation (DRE). Then, a procedure using the analytical strategy for the LQ pursuit-evasion game in arbitrary Keplerian reference orbits is proposed. Simulations validate the effectiveness of the developed analytical strategy, showing that the results obtained from the analytical strategy are in good agreement with those obtained from the numerical one. Furthermore, it is worth emphasizing that in the simulations presented in this paper, the developed analytical strategy saves more than 99.8$\%$ of the CPU time compared to the numerical one, demonstrating significant computational efficiency. The developed strategy is also applicable to the pursuit-evasion game considering orbital disturbances. Finally, comparison between the developed and conventional strategies reveals the advantage of the developed strategy. While the conventional strategy succeeds in only two out of six test cases, the developed strategy achieves success in all six test cases, confirming its effectiveness in pursuit-evasion games in high-eccentricity reference orbits. The analytical strategy can also be applied to other mission scenarios, such as terminal guidance for near-Earth asteroid and space debris defense.

\section*{Declaration of competing interest}
The authors declare that they have no known competing financial interests or personal relationships that could have appeared to influence the work reported in this paper.

\section*{Acknowledgements}

The authors acknowledge the suggestions from Shuo Song and Hongxin Mao from School of Astronautics, Beihang University. The first author acknowledges the financial support from the Outstanding Research Project of Shen Yuan Honors College, BUAA (Grant No. 230122205). The second author acknowledges the financial support from the National Natural Science Foundation of China (Grant No. 12372044). The third author acknowledges the financial support from the National Key Laboratory of Space Intelligent Control (No. HTKJ2024KL502008).

\appendix
\section{Summarization of Acronyms and Variables}\label{secB}

\begin{table}
\caption{Summarization of acronyms.}\label{tabB2}%
\centering
\renewcommand{\arraystretch}{1.5}
\begin{tabular}{@{}ll@{}}
\hline
Acronyms  & Meaning\\
\hline
ARE    & Algebraic Riccati Equation  \\
CW Equations & Clohessy-Wiltshire Equations \\
DRE & Differential Riccati Equation \\
LQ Game & Linear Quadratic Game \\
LVLH Frame & Local-Vertical/Local-Horizontal Frame \\
NEA & Near-Earth Asteroid \\
STM & State Transition Matrix \\
TH Equations & Tschauner-Hempel Equations \\
VOP & Variation of Parameters \\
\hline
\end{tabular}
\end{table}

\begin{table}
\caption{Summarization of variables and definition.}\label{tabB1}%
\centering
\renewcommand{\arraystretch}{1.5}
\begin{tabular}{@{}ll@{}}
\hline
Variable  & Meaning or Definition\\
\hline
$d_c$    & Capture Radius of the LQ Pursuit-Evasion Game  \\
$E$ & Eccentric Anomaly in the Elliptic Reference Orbit \\
$e$ & Eccentricity of the Reference Orbit\\
$f$ & True Anomaly of the Reference Orbit \\
$f_0$ & Initial True Anomaly of the LQ Pursuit-Evasion Game \\
$f_f$ & Terminal True Anomaly of the LQ Pursuit-Evasion Game \\
$\mathcal{H}$ & Hamiltonian Function of the LQ Pursuit-Evasion Game \\
$H$ & Hyperbolic Eccentric Anomaly in the Hyperbolic Reference Orbit \\
$h_f$ & Step Size for Solving the Dynamical Equations \\
$J$ & Cost Function of the LQ Pursuit-Evasion Game \\
$\bm{J}$, $\bm{J}_1$ & Matrices Used in the Calculation of the STMs \\
$n$ & $n=\sqrt{\mu/p^3}$ \\
$\bm{P}$ & Riccati Matrix of the LQ Pursuit-Evasion Game \\
$p$ & Semilatus Rectum of the Reference Orbit \\
$\bm{R}_p$, $\bm{R}_e$, $\bm{S}$ & Weighting Matrices of the LQ Pursuit-Evasion Game \\
$r_p$ & Weight of the Weighting Matrix $\bm{R}_p$ \\
$r_e$ & Weight of the Weighting Matrix $\bm{R}_e$ \\
$\tilde{r}\left(f_f\right)$ & Position Error of the LQ Pursuit-Evasion Game \\
$s_r$ & Weight of the Weighting Matrix $\bm{S}$ \\
$\bm{u}_p$, $\bm{u}_e$ & Control Inputs of the Pursuer and Evader \\
${{\tilde{\bm X}}_p}$ & States of the Pursuer in the LVLH Frame \\
${{\tilde{\bm X}}_e}$ & States of the Evader in the LVLH Frame \\
${{\tilde{\bm x}}}$ & Relative States between the Pursuer and Evader in the LVLH Frame \\
$\beta$ & $\beta  = 1/n^2$\\
$\Delta f$ & Game Time of the LQ Pursuit-Evasion Game \\
$\mu$ & Earth Gravitational Constant \\
$\rho$  & $\rho=1+e\cos{f}$\\
$\Phi$ & Terminal Functional of the LQ Pursuit-Evasion Game \\
$\bm{\Omega}_{11}$, $\bm{\Omega}_{12}$, $\bm{\Omega}_{22}$ & STMs Used in the Derivation of the Analytical Solution of the DRE \\
\hline
\end{tabular}
\end{table}

\section{Components of the Matrix $\bm{J}$ in Elliptic Reference Orbits}\label{secA1}
According to Eq. \eqref{eq27}, the matrix $\bm{J}$ is a symmetric matrix. Therefore, only the upper triangular components are presented. The subscript $\bm{W}_{ij}$ denotes the element of the \textit{i}-th row and \textit{j}-th column of the matrix $\bm{W}$. The components of the matrix $\bm{J}$ calculated in this Appendix are related to Eq. \eqref{eqA} and Section \ref{subsubsec3.1.1}.
\begin{equation}
{\bm{J}_{11}}\left( f \right) =  - \frac{{3 }}{{2{{\left( {1 - {e^2}} \right)}^{\frac{15}{2}}}}}\left( \begin{gathered}
   - \left( {\left( { - \frac{2}{{15}}{e^4} - \frac{2}{{15}}{e^2} + \frac{2}{{15}}} \right)\sin E + {e^3}E} \right){e^3}{\cos ^4}E \hfill \\
   - \frac{8}{3}\left( {\left( { - \frac{5}{{32}}{e^4} - \frac{1}{{16}}{e^2} - \frac{5}{{16}}} \right)\sin E + e\left( {1 + {e^2}} \right)E} \right){e^2}{\cos ^3}E \hfill \\
   - 2e\left( {\left( {\frac{7}{9} - \frac{4}{{45}}{e^6} + \left( {{E^2} - \frac{{44}}{{45}}} \right){e^4} - \frac{7}{{15}}{e^2}} \right)\sin E - \frac{{13}}{2}{e^3}E - 4eE} \right){\cos ^2}E \hfill \\
   + \left( {\left( {1 - \frac{{49}}{{24}}{e^6} + \left( {3{E^2} - \frac{{187}}{{12}}} \right){e^4} - \frac{{85}}{{12}}{e^2}} \right)\sin E - 8e\left( {1 + {e^4}} \right)E} \right)\cos E \hfill \\
   + \left( {\frac{{16}}{{45}}{e^7} + \left( { - 4{E^2} + \frac{{1016}}{{45}}} \right){e^5} + \left( {\frac{{126}}{5} + 6{E^2}} \right){e^3} + \frac{{128}}{9}e} \right)\sin E \hfill \\
   + \left( { - \frac{5}{3} - \frac{{49}}{{24}}{e^6} + \left( {{E^2} - \frac{{251}}{{12}}} \right){e^4} + \left( { - \frac{{63}}{4} - 2{E^2}} \right){e^2}} \right)E \hfill \\ 
\end{gathered}  \right)
\label{eqA1}
\end{equation}

\begin{equation}
{\bm{J}_{12}}(f) = \frac{1}{{60{{\left( {1 - {e^2}} \right)}^6}}}\left( \begin{gathered}
   - 12{e^3}{\cos ^5}E + \left( { - 15{e^6} + 45{e^4} + 75{e^2}} \right){\cos ^4}E + \left( { - 40{e^5} - 100{e^3} - 140e} \right){\cos ^3}E \hfill \\
   + \left( { - 60{e^3}\left( {{e^2} - 2} \right)E\sin E + 150{e^4} + 150{e^2} + 90} \right){\cos ^2}E \hfill \\
   + \left( {\left( {90{e^4}E - 360{e^2}E} \right)\sin E - 120{e^5} + 240{e^3} + 300e} \right)\cos E \hfill \\
   + 45\left( {\left( { - {e^3} + 4e} \right){{\sin }^2}E - \frac{8}{3}\left( {{e^4} - \frac{7}{2}{e^2} - 3} \right)E\sin E + e\left( {{e^2} - 6} \right){E^2}} \right)e \hfill \\ 
\end{gathered}  \right)
\label{eqA2}
\end{equation}

\begin{equation}
{\bm{J}_{13}}(f) = \frac{1}{{60{{\left( {1 - {e^2}} \right)}^6}}}\left( \begin{gathered}
   - 12{e^4}{\cos ^5}E + \left( {30{e^5} + 75{e^3}} \right){\cos ^4}E + \left( { - 100{e^4} - 180{e^2}} \right){\cos ^3}E \hfill \\
   + \left( {60{e^4}E\sin E + 180{e^3} + 210e} \right){\cos ^2}E \hfill \\
   + \left( { - 270{e^3}E\sin E + 120{e^4} + 420{e^2} - 120} \right)\cos E + 135{e^3}{\sin ^2}E \hfill \\
   + \left( {120{e^4}E + 540{e^2}E} \right)\sin E - 135{e^3}{E^2} - 90e{E^2} \hfill \\ 
\end{gathered}  \right)
\label{eqA3}
\end{equation}

\begin{equation}
{\bm{J}_{14}}\left( f \right) =  - \frac{3}{{2{{\left( {1 - {e^2}} \right)}^{\frac{{15}}{2}}}}}\left( \begin{gathered}
  \left( {\left( { - \frac{2}{{15}}{e^4} + \frac{4}{{15}}{e^6}} \right)\sin E - {e^5}E} \right){\cos ^4}E \hfill \\
   + \frac{4}{3}{e^2}{\cos ^3}E\left( {\left( { - \frac{1}{4}{e^5} + \frac{3}{{16}}{e^3} + \frac{9}{8}e} \right)\sin E + E\left( {{e^4} - 4{e^2} - 1} \right)} \right) \hfill \\
   - 2\left( \begin{gathered}
  \left( {\frac{{{e^7}}}{9} - \frac{8}{{45}}{e^5} + {e^3}\left( {{E^2} - \frac{{86}}{{45}}} \right) + \frac{{11}}{9}e} \right)\sin E \hfill \\
   + 2\left( {{e^4} - \frac{{21}}{4}{e^2} - 1} \right)E \hfill \\ 
\end{gathered}  \right)e{\cos ^2}E \hfill \\
   + \left( {\left( {\frac{{11}}{6}{e^7} - \frac{{13}}{8}{e^5} + \left( {3{E^2} - \frac{{275}}{{12}}} \right){e^3} - e} \right)\sin E - 4E{{\left( {1 + {e^2}} \right)}^2}} \right)\cos E \hfill \\
   + \left( {\frac{{16}}{3} - \frac{4}{9}{e^8} - \frac{{298}}{{45}}{e^6} + \left( {\frac{{1304}}{{45}} - 4{E^2}} \right){e^4} + \left( {6{E^2} + \frac{{316}}{9}} \right){e^2}} \right)\sin E \hfill \\
   + eE\left( {\frac{{11}}{6}{e^6} + \frac{{41}}{{24}}{e^4} + \left( {{E^2} - \frac{{129}}{4}} \right){e^2} - 2{E^2} - \frac{{35}}{3}} \right) \hfill \\ 
\end{gathered}  \right)
\label{eqA4}
\end{equation}

\begin{equation}
{\bm{J}_{15}}(f) =0 \text{ }\text{ }\text{ }\text{ }{\bm{J}_{16}}(f) =0
\label{eqA5}
\end{equation}

\begin{equation}
{\bm{J}_{22}}(f) = \frac{{1 }}{{2{{\left( {1 - {e^2}} \right)}^{\frac{11}{2}}}}}\left( \begin{gathered}
  \left( \begin{gathered}
   - \frac{2}{5}{e^3}{\cos ^4}E + \frac{5}{2}{e^2}{\cos ^3}E - \frac{2}{3}e\left( {{e^6} - 5{e^4} + \frac{{39}}{5}{e^2} + 7} \right){\cos ^2}E \hfill \\
   + \left( {{e^6} - 9{e^4} + \frac{{75}}{4}{e^2} + 3} \right)\cos E - \frac{4}{3}{e^7} + \frac{{26}}{3}{e^5} - \frac{{22}}{5}{e^3} - \frac{{82}}{3}e \hfill \\ 
\end{gathered}  \right)\sin E \hfill \\
   + E\left( {{e^6} - 11{e^4} + \frac{{83}}{4}{e^2} + 5} \right) \hfill \\ 
\end{gathered}  \right)\label{eqA6}
\end{equation}

\begin{equation}
{\bm{J}_{23}}\left( f \right) =  - \frac{{3 }}{{2{{\left( {1 - {e^2}} \right)}^{\frac{11}{2}}}}}\left( {\left( \begin{gathered}
  \frac{2}{{15}}{e^4}{\cos ^4}E - \frac{5}{6}{e^3}{\cos ^3}E \hfill \\
   + \left( { - \frac{2}{9}{e^6} + \frac{{28}}{{45}}{e^4} + 2{e^2}} \right){\cos ^2}E \hfill \\
   + \left( {{e^5} - \frac{{13}}{4}{e^3} - \frac{7}{3}e} \right)\cos E \hfill \\
   - \frac{4}{9}{e^6} - \frac{{34}}{{45}}{e^4} + 8{e^2} + \frac{4}{3} \hfill \\ 
\end{gathered}  \right)\sin E + eE\left( {{e^4} - \frac{{31}}{{12}}{e^2} - \frac{{11}}{3}} \right)} \right)\label{eqA7}
\end{equation}

\begin{equation}
{\bm{J}_{24}}\left( f \right) = \frac{1}{{60{{\left( {1 - {e^2}} \right)}^6}}}\left( \begin{gathered}
   - 12{e^4}{\cos ^5}E + \left( { - 30{e^5} + 135{e^3}} \right){\cos ^4}E + \left( {40{e^6} - 100{e^4} - 220{e^2}} \right){\cos ^3}E \hfill \\
   + \left( { - 60{e^2}\left( {{e^2} - 2} \right)E\sin E - 120{e^5} + 465{e^3} - 90e} \right){\cos ^2}E \hfill \\
   + \left( {\left( {90{e^3}E - 360eE} \right)\sin E - 60{e^2} + 480} \right)\cos E - 120\left( {{e^4} - \frac{7}{2}{e^2} - 3} \right)E\sin E \hfill \\
   + 45\left( {\left( {{E^2} - 1} \right){e^2} - 6{E^2} + 4} \right)e \hfill \\ 
\end{gathered}  \right)\label{eqA8}
\end{equation}

\begin{equation}
{\bm{J}_{25}}(f) =0 \text{ }\text{ }\text{ }\text{ }{\bm{J}_{26}}(f) =0
\label{eqA9}
\end{equation}

\begin{equation}
{\bm{J}_{33}}\left( f \right) = \frac{1}{{120{{\left( {1 - {e^2}} \right)}^{\frac{{11}}{2}}}}}\left( \begin{gathered}
  \left( \begin{gathered}
   - 24{e^5}{\cos ^4}E + 150{e^4}{\cos ^3}E + \left( { - 32{e^5} - 400{e^3}} \right){\cos ^2}E \hfill \\
   + \left( {225{e^4} + 600{e^2}} \right)\cos E - 64{e^5} - 800{e^3} - 600e \hfill \\ 
\end{gathered}  \right)\sin E \hfill \\
   + 225{e^4}E + 600{e^2}E + 120E \hfill \\ 
\end{gathered}  \right)\label{eqA10}
\end{equation}

\begin{equation}
{\bm{J}_{34}}\left( f \right) = \frac{1}{{60{{\left( {1 - {e^2}} \right)}^6}}}\left( \begin{gathered}
   - 12{e^5}{\cos ^5}E + \left( { - 15{e^6} + 120{e^4}} \right){\cos ^4}E + \left( {60{e^5} - 340{e^3}} \right)\cos 3E \hfill \\
   + \left( {60{e^3}E\sin E - 90{e^4} + 345{e^2}} \right){\cos ^2}E + \left( { - 270{e^2}E{{\sin }^2}E + 180{e^3} + 240e} \right)\cos E \hfill \\
   + \left( {120{e^3}E + 540eE} \right)\sin E + \left( { - 135{E^2} + 135} \right){e^2} - 90{E^2} \hfill \\ 
\end{gathered}  \right)\label{eqA11}
\end{equation}

\begin{equation}
{\bm{J}_{35}}(f) =0 \text{ }\text{ }\text{ }\text{ }{\bm{J}_{36}}(f) =0
\label{eqA12}
\end{equation}

\begin{equation}
{\bm{J}_{44}}\left( f \right) = \frac{{3eE}}{{{{\left( {1 - {e^2}} \right)}^{\frac{{15}}{2}}}}}\left( \begin{gathered}
  \left( \begin{gathered}
   - \frac{1}{{15}}{e^6}{\cos ^4}E + \left( { - \frac{1}{6}{e^7} + \frac{5}{4}{e^5} - \frac{{43}}{{24}}{e^3}} \right){\cos ^3}E \hfill \\
   + {e^2}\left( { - \frac{1}{9}{e^6} + \frac{{17}}{{15}}{e^4} + {E^2} - \frac{{28}}{9}{e^2} + \frac{4}{3}} \right){\cos ^2}E \hfill \\
   - \frac{3}{2}e\left( { - \frac{7}{{18}}{e^6} + \frac{{37}}{{12}}{e^4} + {E^2} - \frac{{287}}{{72}}{e^2} - \frac{{119}}{{18}}} \right)\cos E \hfill \\
   - \frac{2}{9}{e^8} - \frac{{{e^6}}}{{15}} + \frac{{103}}{9}{e^4} + \left( {2{E^2} - \frac{{79}}{3}} \right){e^2} - 3{E^2} - 16 \hfill \\ 
\end{gathered}  \right)\sin E \hfill \\
   + \left( \begin{gathered}
   - {e^4}{\cos ^4}E + \left( {\frac{8}{3}{e^5} - 8{e^3}} \right){\cos ^3}E + \left( { - 8{e^4} + 29{e^2}} \right){\cos ^2}E \hfill \\
   - 16e\cos E - \frac{7}{6}{e^8} + \frac{{29}}{4}{e^6} + \frac{{27}}{8}{e^4} + \left( { - \frac{{283}}{6} + {E^2}} \right){e^2} - \frac{8}{3} - 2{E^2} \hfill \\ 
\end{gathered}  \right) \hfill \\ 
\end{gathered}  \right)\label{eqA13}
\end{equation}

\begin{equation}
{\bm{J}_{45}}(f) =0 \text{ }\text{ }\text{ }\text{ }{\bm{J}_{46}}(f) =0
\label{eqA14}
\end{equation}

\begin{equation}
{\bm{J}_{55}}\left( f \right) =  - \frac{3}{{8{{\left( {1 - {e^2}} \right)}^{\frac{9}{2}}}}}\left( {\left( \begin{gathered}
  \frac{8}{{15}}{e^3}{\cos ^4}E - 2{e^2}{\cos ^3}E + \left( { - \frac{8}{{45}}{e^3} + \frac{8}{3}e} \right){\cos ^2}E \hfill \\
   + \left( {{e^2} - \frac{4}{3}} \right)\cos E - \frac{{16}}{{45}}{e^3} - \frac{8}{3}e \hfill \\ 
\end{gathered}  \right)\sin E + E\left( {{e^2} + \frac{4}{3}} \right)} \right)\label{eqA15}
\end{equation}

\begin{equation}
{\bm{J}_{56}}\left( f \right) =  - \frac{{\cos E}}{{4{{\left( {1 - {e^2}} \right)}^5}}}\left( \begin{gathered}
   - \frac{4}{5}{e^3}{\cos ^4}E + \left( {{e^4} + 3{e^2}} \right){\cos ^3}E \hfill \\
   + \left( { - 4{e^3} - 4e} \right){\cos ^2}E + \left( {6{e^2} + 2} \right)\cos E - 4e \hfill \\ 
\end{gathered}  \right)\label{eqA16}
\end{equation}

\begin{equation}
{\bm{J}_{66}}\left( f \right) = \frac{9}{{4{{\left( {1 - {e^2}} \right)}^{\frac{{11}}{2}}}}}\left( \begin{gathered}
  \left( \begin{gathered}
   - \frac{4}{{45}}{e^3}{\cos ^4}E + \left( {\frac{2}{9}{e^4} + \frac{1}{3}{e^2}} \right){\cos ^3}E + \left( { - \frac{4}{{27}}{e^5} - \frac{{136}}{{135}}{e^3} - \frac{4}{9}e} \right){\cos ^2}E \hfill \\
   + \left( {{e^4} + \frac{{11}}{6}{e^2} + \frac{2}{9}} \right)\cos E - \frac{8}{{27}}{e^5} - \frac{{452}}{{135}}{e^3} - \frac{{16}}{9}e \hfill \\ 
\end{gathered}  \right)\sin E \hfill \\
   + E\left( {{e^4} + \frac{{41}}{{18}}{e^2} + \frac{2}{9}} \right) \hfill \\ 
\end{gathered}  \right)\label{eqA17}
\end{equation}

\section{Components of the Matrix $\bm{J}$ in Parabolic Reference Orbits}\label{secA2}
The components of the matrix $\bm{J}$ calculated in this Appendix are related to Section \ref{subsubsec3.1.2}.
\begin{equation}
{\bm{J}_{11}}\left( f \right) =  - \frac{{\sin \left( {\frac{f}{2}} \right)}}{{144144000{{\cos }^{15}}\left( {\frac{f}{2}} \right)}}\left( \begin{gathered}
  237568{\cos ^{14}}\left( {\frac{f}{2}} \right) + 118784{\cos ^{12}}\left( {\frac{f}{2}} \right) + 89088{\cos ^{10}}\left( {\frac{f}{2}} \right) \hfill \\
   - 1847680{\cos ^8}\left( {\frac{f}{2}} \right) + 185080{\cos ^6}\left( {\frac{f}{2}} \right) + 166572{\cos ^4}\left( {\frac{f}{2}} \right) \hfill \\
   - 72534{\cos ^2}\left( {\frac{f}{2}} \right) - 3003 \hfill \\ 
\end{gathered}  \right)\label{eqB1}
\end{equation}

\begin{equation}
{\bm{J}_{12}}\left( f \right) = \frac{1}{{25{{\left( {1 + \cos f} \right)}^5}}} + \frac{3}{{20{{\left( {1 + \cos f} \right)}^4}}} - \frac{4}{{15{{\left( {1 + \cos f} \right)}^3}}} - \frac{1}{{15{{\left( {1 + \cos f} \right)}^6}}}\label{eqB2}
\end{equation}

\begin{equation}
{\bm{J}_{13}}\left( f \right) = \frac{1}{{10{{\left( {1 + \cos f} \right)}^4}}} + \frac{2}{{25{{\left( {1 + \cos f} \right)}^5}}} - \frac{1}{{15{{\left( {1 + \cos f} \right)}^6}}}\label{eqB3}
\end{equation}

\begin{equation}
{\bm{J}_{14}}\left( f \right) =  - \frac{{\sin \left( {\frac{f}{2}} \right)}}{{48048000{{\cos }^{15}}\left( {\frac{f}{2}} \right)}}\left( \begin{gathered}
  51456{\cos ^{14}}\left( {\frac{f}{2}} \right) + 25728{\cos ^{12}}\left( {\frac{f}{2}} \right) + 19296{\cos ^{10}}\left( {\frac{f}{2}} \right) \hfill \\
   + 977040{\cos ^8}\left( {\frac{f}{2}} \right) + 104160{\cos ^6}\left( {\frac{f}{2}} \right) - 26376{\cos ^4}\left( {\frac{f}{2}} \right) \hfill \\
   - 24178{\cos ^2}\left( {\frac{f}{2}} \right) - 1001 \hfill \\ 
\end{gathered}  \right)\label{eqB4}
\end{equation}

\begin{equation}
{\bm{J}_{15}}(f) =0 \text{ }\text{ }\text{ }\text{ }{\bm{J}_{16}}(f) =0
\label{eqB5}
\end{equation}

\begin{equation}
{\bm{J}_{22}}\left( f \right) = \frac{1}{{352}}{\tan ^{11}}\left( {\frac{f}{2}} \right) + \frac{5}{{288}}{\tan ^9}\left( {\frac{f}{2}} \right) + \frac{1}{{112}}{\tan ^7}\left( {\frac{f}{2}} \right) + \frac{1}{{80}}{\tan ^5}\left( {\frac{f}{2}} \right) + \frac{{13}}{{96}}{\tan ^3}\left( {\frac{f}{2}} \right) + \frac{9}{{32}}\tan \left( {\frac{f}{2}} \right)\label{eqB6}
\end{equation}

\begin{equation}
{\bm{J}_{23}}\left( f \right) = \frac{1}{{352}}{\tan ^{11}}\left( {\frac{f}{2}} \right) + \frac{5}{{288}}{\tan ^9}\left( {\frac{f}{2}} \right) + \frac{3}{{112}}{\tan ^7}\left( {\frac{f}{2}} \right) - \frac{1}{{80}}{\tan ^5}\left( {\frac{f}{2}} \right) - \frac{7}{{96}}{\tan ^3}\left( {\frac{f}{2}} \right) - \frac{3}{{32}}\tan \left( {\frac{f}{2}} \right)\label{eqB7}
\end{equation}

\begin{equation}
{\bm{J}_{24}}\left( f \right) =  - \frac{4}{{25{{\left( {1 + \cos f} \right)}^5}}} - \frac{1}{{10{{\left( {1 + \cos f} \right)}^4}}} + \frac{2}{{5{{\left( {1 + \cos f} \right)}^3}}} - \frac{1}{{15{{\left( {1 + \cos f} \right)}^6}}}\label{eqB8}
\end{equation}

\begin{equation}
{\bm{J}_{25}}(f) =0 \text{ }\text{ }\text{ }\text{ }{\bm{J}_{26}}(f) =0
\label{eqB9}
\end{equation}

\begin{equation}
{\bm{J}_{33}}\left( f \right) = \frac{1}{{352}}{\tan ^{11}}\left( {\frac{f}{2}} \right) + \frac{5}{{288}}{\tan ^9}\left( {\frac{f}{2}} \right) + \frac{5}{{112}}{\tan ^7}\left( {\frac{f}{2}} \right) + \frac{1}{{16}}{\tan ^5}\left( {\frac{f}{2}} \right) + \frac{5}{{96}}{\tan ^3}\left( {\frac{f}{2}} \right) + \frac{1}{{32}}\tan \left( {\frac{f}{2}} \right)\label{eqB10}
\end{equation}

\begin{equation}
{\bm{J}_{34}}\left( f \right) =  - \frac{3}{{25{{\left( {1 + \cos f} \right)}^5}}} - \frac{3}{{20{{\left( {1 + \cos f} \right)}^4}}} -  - \frac{1}{{15{{\left( {1 + \cos f} \right)}^6}}}\label{eqB11}
\end{equation}

\begin{equation}
{\bm{J}_{35}}(f) =0 \text{ }\text{ }\text{ }\text{ }{\bm{J}_{36}}(f) =0
\label{eqB12}
\end{equation}

\begin{equation}
{\bm{J}_{44}}\left( f \right) = \frac{{\sin \left( {\frac{f}{2}} \right)}}{{144144000{{\cos }^{15}}\left( {\frac{f}{2}} \right)}}\left( \begin{gathered}
  1759232{\cos ^{14}}\left( {\frac{f}{2}} \right) + 879616{\cos ^{12}}\left( {\frac{f}{2}} \right) + 659712{\cos ^{10}}\left( {\frac{f}{2}} \right) \hfill \\
   + 4874080{\cos ^8}\left( {\frac{f}{2}} \right) + 1562120{\cos ^6}\left( {\frac{f}{2}} \right) + 324828{\cos ^4}\left( {\frac{f}{2}} \right) \hfill \\
   + 72534{\cos ^2}\left( {\frac{f}{2}} \right) + 3003 \hfill \\ 
\end{gathered}  \right)\label{eqB13}
\end{equation}

\begin{equation}
{\bm{J}_{45}}(f) =0 \text{ }\text{ }\text{ }\text{ }{\bm{J}_{46}}(f) =0
\label{eqB14}
\end{equation}

\begin{equation}
{\bm{J}_{55}}\left( f \right) = \frac{1}{{72}}{\tan ^9}\left( {\frac{f}{2}} \right) + \frac{3}{{56}}{\tan ^7}\left( {\frac{f}{2}} \right) + \frac{3}{{40}}{\tan ^5}\left( {\frac{f}{2}} \right) + \frac{1}{{24}}{\tan ^3}\left( {\frac{f}{2}} \right)\label{eqB15}
\end{equation}

\begin{equation}
{\bm{J}_{56}}\left( f \right) = \frac{1}{{5{{\left( {1 + \cos f} \right)}^5}}} - \frac{1}{{4{{\left( {1 + \cos f} \right)}^4}}}\label{eqB16}
\end{equation}

\begin{equation}
{\bm{J}_{66}}\left( f \right) = \frac{1}{{352}}{\tan ^{11}}\left( {\frac{f}{2}} \right) + \frac{1}{{288}}{\tan ^9}\left( {\frac{f}{2}} \right) - \frac{1}{{112}}{\tan ^7}\left( {\frac{f}{2}} \right) - \frac{1}{{80}}{\tan ^5}\left( {\frac{f}{2}} \right) + \frac{1}{{96}}{\tan ^3}\left( {\frac{f}{2}} \right) + \frac{1}{{32}}\tan \left( {\frac{f}{2}} \right)\label{eqB17}
\end{equation}

\section{Components of the Matrix $\bm{J}$ in Hyperbolic Reference Orbits}\label{secA3}
The components of the matrix $\bm{J}$ calculated in this Appendix are related to Eq. \eqref{eqB} and Section \ref{subsubsec3.1.3}.
\begin{equation}
{\bm{J}_{11}}\left( f \right) =  - \frac{3}{{2{{\left( {{e^2} - 1} \right)}^{\frac{{15}}{2}}}}}\left( \begin{gathered}
  {e^3}\left( {\left( { - \frac{2}{{15}}{e^4} - \frac{2}{{15}}{e^2} + \frac{2}{{15}}} \right)\sinh H + {e^3}H} \right){\cosh ^4}H \hfill \\
   + \frac{8}{3}{e^2}{\cosh ^3}H\left( {\left( { - \frac{5}{{32}}{e^4} - \frac{1}{{16}}{e^2} - \frac{5}{{16}}} \right)\sinh H + e\left( {{e^2} + 1} \right)H} \right) \hfill \\
   - 2e\left( {\left( { - \frac{7}{9} + \frac{4}{{45}}{e^6} + \left( {{H^2} + \frac{{44}}{{45}}} \right){e^4} + \frac{7}{{15}}{e^2}} \right)\sinh H + \frac{{13}}{2}{e^3}H + 4eH} \right){\cosh ^2}H \hfill \\
   + \left( {\left( { - 1 + \frac{{49}}{{24}}{e^6} + \left( {3{H^2} + \frac{{187}}{{12}}} \right){e^4} + \frac{{85}}{{12}}{e^2}} \right)\sinh H + 8e\left( {{e^4} + 1} \right)H} \right)\cosh H \hfill \\
   + \left( { - \frac{{16}}{{45}}{e^7} + \left( { - 4{H^2} - \frac{{1016}}{{45}}} \right){e^5} + \left( {6{H^2} - \frac{{126}}{5}} \right){e^3} - \frac{{128}}{9}e} \right)\sinh H \hfill \\
   + H\left( {\frac{5}{3} + \frac{{49}}{{24}}{e^6} + \left( {{H^2} + \frac{{251}}{{12}}} \right){e^4} + \left( { - 2{H^2} + \frac{{63}}{4}} \right){e^2}} \right) \hfill \\ 
\end{gathered}  \right)\label{eqC1}
\end{equation}

\begin{equation}
{\bm{J}_{12}}\left( f \right) = \frac{1}{{60{{\left( {{e^2} - 1} \right)}^6}}}\left( \begin{gathered}
   - 12{e^3}{\cosh ^5}H + \left( { - 15{e^6} + 45{e^4} + 75{e^2}} \right){\cosh ^4}H + \left( { - 40{e^5} - 100{e^3} - 140e} \right){\cosh ^3}H \hfill \\
   + \left( {60{e^3}\left( {{e^2} - 2} \right)H\sinh H + 195{e^4} - 30{e^2} + 90} \right){\cosh ^2}H \hfill \\
   + \left( {\left( { - 90{e^4}H + 360{e^2}H} \right)\sinh H - 120{e^5} + 240{e^3} + 300e} \right)\cosh H \hfill \\
   - 45eH\left( {\left( { - \frac{8}{3}{e^4} + \frac{{28}}{3}{e^2} + 8} \right)\sinh H + e\left( {{e^2} - 6} \right)H} \right) \hfill \\ 
\end{gathered}  \right)\label{eqC2}
\end{equation}

\begin{equation}
{\bm{J}_{13}}\left( f \right) = \frac{1}{{60{{\left( {{e^2} - 1} \right)}^6}}}\left( \begin{gathered}
   - 12{e^4}{\cosh ^5}H + \left( {30{e^5} + 75{e^3}} \right){\cosh ^4}H + \left( { - 100{e^4} - 180{e^2}} \right){\cosh ^3}H \hfill \\
   + \left( { - 60{e^4}H\sinh H + 45{e^3} + 210e} \right){\cosh ^2}H \hfill \\
   + \left( {270{e^3}H\sinh H + 120{e^4} + 420{e^2} - 120} \right)\cosh H \hfill \\
   + 135eH\left( {\left( { - \frac{8}{9}{e^3} - 4e} \right)\sinh H + H\left( {{e^2} + \frac{2}{3}} \right)} \right) \hfill \\ 
\end{gathered}  \right)\label{eqC3}
\end{equation}

\begin{equation}
{\bm{J}_{14}}\left( f \right) =  - \frac{3}{{2{{\left( {{e^2} - 1} \right)}^{\frac{{15}}{2}}}}}\left( \begin{gathered}
  \left( {\left( { - \frac{4}{{15}}{e^6} + \frac{2}{{15}}{e^4}} \right)\sinh H + {e^5}H} \right){\cosh ^4}H \hfill \\
   - \frac{4}{3}{e^2}{\cosh ^3}H\left( {\left( { - \frac{1}{4}{e^5} + \frac{3}{{16}}{e^3} + \frac{9}{8}e} \right)\sinh H + H\left( {{e^4} - 4{e^2} - 1} \right)} \right) \hfill \\
   - 2\left( {\left( { - \frac{{{e^7}}}{9} + \frac{8}{{45}}{e^5} + \left( {{H^2} + \frac{{86}}{{45}}} \right){e^3} - \frac{{11}}{9}e} \right)\sinh H - 2H\left( {{e^4} - \frac{{21}}{4}{e^2} - 1} \right)} \right)e{\cosh ^2}H \hfill \\
   + \left( {\left( { - \frac{{11}}{6}{e^7} + \frac{{13}}{8}{e^5} + \left( {\frac{{275}}{{12}} + 3{H^2}} \right){e^3} + e} \right)\sinh H + 4{{\left( {{e^2} + 1} \right)}^2}H} \right)\cosh H \hfill \\
   + \left( { - \frac{{16}}{3} + \frac{4}{9}{e^8} + \frac{{298}}{{45}}{e^6} + \left( { - 4{H^2} - \frac{{1304}}{{45}}} \right){e^4} + \left( { - \frac{{316}}{9} + 6{H^2}} \right){e^2}} \right)\sinh H \hfill \\
   + eH\left( { - \frac{{11}}{6}{e^6} - \frac{{41}}{{24}}{e^4} + \left( {\frac{{129}}{4} + {H^2}} \right){e^2} + \frac{{35}}{3} - 2{H^2}} \right) \hfill \\ 
\end{gathered}  \right)\label{eqC4}
\end{equation}

\begin{equation}
{\bm{J}_{15}}(f) =0 \text{ }\text{ }\text{ }\text{ }{\bm{J}_{16}}(f) =0
\label{eqC5}
\end{equation}

\begin{equation}
{\bm{J}_{22}}\left( f \right) =  - \frac{1}{{2{{\left( {{e^2} - 1} \right)}^{\frac{{11}}{2}}}}}\left( \begin{gathered}
  \left( \begin{gathered}
   - \frac{2}{5}{e^3}{\cosh ^4}H + \frac{5}{2}{e^2}{\cosh ^3}H - \frac{2}{3}e{\cosh ^2}H\left( {{e^6} - 5{e^4} + \frac{{39}}{5}{e^2} + 7} \right) \hfill \\
   + \left( {{e^6} - 9{e^4} + \frac{{75}}{4}{e^2} + 3} \right)\cosh H - \frac{4}{3}{e^7} + \frac{{26}}{3}{e^5} - \frac{{22}}{5}{e^3} - \frac{{82}}{3}e \hfill \\ 
\end{gathered}  \right)\sinh H \hfill \\
   + H\left( {{e^6} - 11{e^4} + \frac{{83}}{4}{e^2} + 5} \right) \hfill \\ 
\end{gathered}  \right)\label{eqC6}
\end{equation}

\begin{equation}
{\bm{J}_{23}}\left( f \right) = \frac{3}{{2{{\left( {{e^2} - 1} \right)}^{\frac{{11}}{2}}}}}\left( {\left( \begin{gathered}
  \frac{2}{{15}}{e^4}{\cosh ^4}H - \frac{5}{6}{e^3}{\cosh ^3}H \hfill \\
   + \left( { - \frac{2}{9}{e^6} + \frac{{28}}{{45}}{e^4} + 2{e^2}} \right){\cosh ^2}H \hfill \\
   + \left( {{e^5} - \frac{{13}}{4}{e^3} - \frac{7}{3}e} \right)\cosh H \hfill \\
   - \frac{4}{9}{e^6} - \frac{{34}}{{45}}{e^4} + 8{e^2} + \frac{4}{3} \hfill \\ 
\end{gathered}  \right)\sinh H + eH\left( {{e^4} - \frac{{31}}{{12}}{e^2} - \frac{{11}}{3}} \right)} \right)
\label{eqC7}
\end{equation}

\begin{equation}
{\bm{J}_{24}}\left( f \right) = \frac{1}{{60{{\left( {{e^2} - 1} \right)}^6}}}\left( \begin{gathered}
   - 12{e^4}{\cosh ^5}H + \left( { - 30{e^5} + 135{e^3}} \right){\cosh ^4}H + \left( {40{e^6} - 100{e^4} - 220{e^2}} \right){\cosh ^3}H \hfill \\
   + \left( {60{e^2}\left( {{e^2} - 2} \right)H\sinh H - 120{e^5} + 465{e^3} - 90e} \right){\cosh ^2}H \hfill \\
   + \left( {\left( { - 90{e^3}H + 360eH} \right)\sinh H - 60{e^2} + 480} \right)\cosh H \hfill \\
   - 45H\left( {\left( { - \frac{8}{3}{e^4} + \frac{{28}}{3}{e^2} + 8} \right)\sinh H + e\left( {{e^2} - 6} \right)H} \right) \hfill \\ 
\end{gathered}  \right)
\label{eqC8}
\end{equation}

\begin{equation}
{\bm{J}_{25}}(f) =0 \text{ }\text{ }\text{ }\text{ }{\bm{J}_{26}}(f) =0
\label{eqC9}
\end{equation}

\begin{equation}
{\bm{J}_{33}}\left( f \right) = \frac{1}{{120{{\left( {{e^2} - 1} \right)}^{\frac{{11}}{2}}}}}\left( \begin{gathered}
  \left( \begin{gathered}
  24{e^5}{\cosh ^4}H - 150{e^4}{\cosh ^3}H + \left( {32{e^5} + 400{e^3}} \right){\cosh ^2}H \hfill \\
   + \left( { - 225{e^4} - 600{e^2}} \right)\cosh H + 64{e^5} + 800{e^3} + 600e \hfill \\ 
\end{gathered}  \right)\sinh H \hfill \\
   - 225{e^4}H - 600{e^2}H - 120H \hfill \\ 
\end{gathered}  \right)\label{eqC10}
\end{equation}

\begin{equation}
{\bm{J}_{34}}\left( f \right) = \frac{1}{{60{{\left( {{e^2} - 1} \right)}^6}}}\left( \begin{gathered}
   - 12{e^5}{\cosh ^5}H + \left( { - 15{e^6} + 120{e^4}} \right){\cosh ^4}H + \left( {60{e^5} - 340{e^3}} \right){\cosh ^3}H \hfill \\
   + \left( { - 60{e^3}H\sinh H - 90{e^4} + 345{e^2}} \right){\cosh ^2}H \hfill \\
   + \left( {270{e^2}H\sinh H + 180{e^3} + 240e} \right)\cosh H \hfill \\
   + 135H\left( {\left( { - \frac{8}{9}{e^3} - 4e} \right)\sinh H + H\left( {{e^2} + \frac{2}{3}} \right)} \right) \hfill \\ 
\end{gathered}  \right)\label{eqC11}
\end{equation}

\begin{equation}
{\bm{J}_{35}}(f) =0 \text{ }\text{ }\text{ }\text{ }{\bm{J}_{36}}(f) =0
\label{eqC12}
\end{equation}

\begin{equation}
{\bm{J}_{44}}\left( f \right) = \frac{{3e}}{{{{\left( {{e^2} - 1} \right)}^{\frac{{15}}{2}}}}}\left( \begin{gathered}
  \left( \begin{gathered}
  \frac{1}{{15}}{e^6}{\cosh ^4}H + \left( {\frac{1}{6}{e^7} - \frac{5}{4}{e^5} + \frac{{43}}{{24}}{e^3}} \right){\cosh ^3}H \hfill \\
   + {e^2}\left( {\frac{1}{9}{e^6} - \frac{{17}}{{15}}{e^4} + {H^2} + \frac{{28}}{9}{e^2} - \frac{4}{3}} \right){\cosh ^2}H \hfill \\
   - \frac{3}{2}e\cosh H\left( {\frac{7}{{18}}{e^6} - \frac{{37}}{{12}}{e^4} + {H^2} + \frac{{287}}{{72}}{e^2} + \frac{{119}}{{18}}} \right)e\cosh H \hfill \\
   + \frac{2}{9}{e^8} + \frac{{{e^6}}}{{15}} - \frac{{103}}{9}{e^4} + \left( {2{H^2} + \frac{{79}}{3}} \right){e^2} - 3{H^2} + 16 \hfill \\ 
\end{gathered}  \right)\sinh H \hfill \\
   + H\left( \begin{gathered}
  {e^4}{\cosh ^4}H + \left( {8{e^3} - \frac{8}{3}{e^5}} \right){\cosh ^3}H + \left( {8{e^4} - 29{e^2}} \right){\cosh ^2}H \hfill \\
   + 16e\cosh H + \frac{7}{6}{e^8} - \frac{{29}}{4}{e^6} - \frac{{27}}{8}{e^4} + \left( {{H^2} + \frac{{283}}{6}} \right){e^2} + \frac{8}{3} - 2{H^2} \hfill \\ 
\end{gathered}  \right) \hfill \\ 
\end{gathered}  \right)\label{eqC13}
\end{equation}

\begin{equation}
{\bm{J}_{45}}(f) =0 \text{ }\text{ }\text{ }\text{ }{\bm{J}_{46}}(f) =0
\label{eqC14}
\end{equation}

\begin{equation}
{\bm{J}_{55}}\left( f \right) = \frac{3}{{8{{\left( {{e^2} - 1} \right)}^{\frac{9}{2}}}}}\left( {\left( \begin{gathered}
  \frac{8}{{15}}{e^3}{\cosh ^4}H - 2{e^2}{\cosh ^3}H + \left( { - \frac{8}{{45}}{e^3} + \frac{8}{3}e} \right){\cosh ^2}H \hfill \\
   + \left( {{e^2} - \frac{4}{3}} \right)\cosh H - \frac{{16}}{{45}}{e^3} - \frac{8}{3}e \hfill \\ 
\end{gathered}  \right)\sinh H + H\left( {{e^2} + \frac{4}{3}} \right)} \right)\label{eqC15}
\end{equation}

\begin{equation}
{\bm{J}_{56}}\left( f \right) =  - \frac{{\cosh H}}{{4{{\left( {{e^2} - 1} \right)}^5}}}\left( \begin{gathered}
   - \frac{4}{5}{e^3}{\cosh ^4}H + \left( {{e^4} + 3{e^2}} \right){\cosh ^3}H \hfill \\
   + \left( { - 4{e^3} - 4e} \right){\cosh ^2}H + \left( {6{e^2} + 2} \right)\cosh H - 4e \hfill \\ 
\end{gathered}  \right)\label{eqC16}
\end{equation}

\begin{equation}
{\bm{J}_{66}}\left( f \right) =  - \frac{9}{{4{{\left( {{e^2} - 1} \right)}^{\frac{{11}}{2}}}}}\left( {\left( \begin{gathered}
   - \frac{4}{{45}}{e^3}{\cosh ^4}H + \left( {\frac{2}{9}{e^4} + \frac{1}{3}{e^2}} \right){\cosh ^3}H \hfill \\
   + \left( { - \frac{4}{{27}}{e^5} - \frac{{136}}{{135}}{e^3} - \frac{4}{9}e} \right){\cosh ^2}H \hfill \\
   + \left( {{e^4} + \frac{{11}}{6}{e^2} + \frac{2}{9}} \right)\cosh H - \frac{8}{{27}}{e^5} - \frac{{452}}{{135}}{e^3} - \frac{{16}}{9}e \hfill \\ 
\end{gathered}  \right)\sinh H + H\left( {{e^4} + \frac{{41}}{{18}}{e^2} + \frac{2}{9}} \right)} \right)\label{eqC17}
\end{equation}

\section{Analytical Strategy for Defense of Near-Earth Asteroids and Space Debris}\label{secA4}

To further illustrate the applicability of the analytical strategy developed in this paper, the scenario is extended to the defense of near-Earth asteroids or space debris. In this extension, the pursuer denotes the spacecraft and the evader denotes the near-Earth asteroid or space debris. Since the near-Earth asteroids or space debris may lack maneuvering capabilities, the dynamical equations presented in Eq. \eqref{eq7} require modification to account for the absence of evader control inputs:

\begin{equation}
\left\{ \begin{gathered}
  {{\tilde{\bm X}'}_p} = \bm{A}{{\tilde{\bm X}}_p} + \bm{B}{\bm{u}_p} \hfill \\
  {{\tilde{\bm X}'}_e} = \bm{A}{{\tilde{\bm X}}_e} \hfill \\ 
\end{gathered}  \right.\label{eqD1}
\end{equation}
Then, the dynamical equations between the pursuer and evader can be rewritten as:

\begin{equation}
\tilde{\bm x}' = \bm{A}\tilde{\bm x} + \bm{B}{\bm{u}_p} \label{eqD2}
\end{equation}
The cost function and Hamiltonian can be expressed as:

\begin{equation}
J = \frac{1}{2}\left( {\tilde{\bm x}{{\left( {{f_f}} \right)}^{\text{T}}}\bm{S}\tilde{\bm x}\left( {{f_f}} \right)} \right) + \frac{1}{2}\int_{f_0}^{{f_f}} { {{\bm{u}_p}^{\text{T}}{\bm{R}_p}{\bm{u}_p} } {\text{d}}f} \label{eqD3}
\end{equation}

\begin{equation}
\mathcal{H} = \frac{1}{2} {{\bm{u}_p}^{\text{T}}{\bm{R}_p}{\bm{u}_p}} + {\bm{\lambda }^{\text{T}}}\left( {\bm{A}\tilde{\bm x} + \bm{B}{\bm{u}_p} } \right)\label{eqD4}
\end{equation}

With the same terminal condition presented in Eq. \eqref{eq19}, the optimal strategy of the spacecraft is expressed as:

\begin{equation}
{\bm{u}_p} =  - {\left( {{\bm{R}_p}} \right)^{ - 1}}{\bm{B}^{\text{T}}}\bm{\lambda }\label{eqD5}
\end{equation}
Based on Eq. \eqref{eqD5}, the DRE is written as:
\cite{liao2021research}:
\begin{equation}
\bm{P}'\left( f \right)\bm{ = } - {\bm{A}^{\text{T}}}\bm{P}\left( f \right) - \bm{P}\left( f \right)\bm{A} + \bm{P}\left( f \right)\bm{B} {{{\left( {{\bm{R}_p}} \right)}^{ - 1}}} {\bm{B}^{\text{T}}}\bm{P}\left( f \right)\label{eqD6}
\end{equation}

\begin{equation}
\bm{P}\left( {{f_f}} \right) = \bm{S}\label{eqD7}
\end{equation}
Then, the matrix ${\bm{\Omega }_{12}}\left( {{f_2},{\text{ }}{f_1}} \right)$ is expressed as:

\begin{equation}
{\bm{\Omega }_{12}}\left( {{f_2},{\text{ }}{f_1}} \right) = -\frac{1}{{{n^4}}} { \frac{1}{{{r_p}}}} \bm{J}_1 \label{eqD70}
\end{equation}
Therefore, the new matrix $\bm{J}_1$ for the extended scenario can be obtained using the same method proposed in this paper. Therefore, based on aforementioned discussion, the analytical strategy for the defense of near-Earth asteroids or space debris can be also obtained.



\printcredits

\end{document}